\documentclass[12pt]{iopart}

\usepackage{graphicx}
\usepackage{amssymb}
\usepackage{bm}

\expandafter\let\csname equation*\endcsname\relax

\expandafter\let\csname endequation*\endcsname\relax

\expandafter\let\csname array*\endcsname\relax

\expandafter\let\csname endarray*\endcsname\relax

\usepackage{amsmath}
\usepackage{color}
\usepackage{subfigure}
\usepackage{appendix}

\numberwithin{equation}{section}

\def\bea{\begin{eqnarray}}
\def\eea{\end{eqnarray}}
\def\ee{\mbox{e}}

\def\sech{\mbox{sech}}

\def\dd{\mbox{d}}
\def\eps{\epsilon}

\newcommand{\cmmnt}[1]{\ignorespaces}
\def\ri{\mathrm{i}}

\def\ye{Y}
\def\Ue{U}
\def\hy{y}

\def\asm{\mu}

\def\dg{\delta}
\def\sq{\sqrt{2}}

\def\tnh{\tanh y}
\def\ssq{\mbox{sech}^2 y}

\begin{document}

\title[Paper 1]{Termination points and homoclinic glueing for a class of inhomogeneous nonlinear ordinary differential equations}

\author{J. S. Keeler$^1$, M. G. Blyth$^2$, J. R. King$^3$}

\address{$^1$ School of Mathematics, University of Manchester, Alan Turing Building, Oxford Road, Manchester, M13 9PL, UK, $^2$ School of Mathematics, University of East Anglia, Norwich Research Park, Norwich, NR4 7TJ, UK, $^3$ School of Mathematical Sciences, University of Nottingham, Nottingham, NG7 2RD, UK}
\ead{jack.keeler@manchester.ac.uk, m.blyth@uea.ac.uk, John.King@uea.ac.uk}
\vspace{10pt}
\begin{indented}
\item[]August 2019
\end{indented}

\begin{abstract}
Solutions $u(x)$ to the class of inhomogeneous nonlinear ordinary differential equations taking the form
\bea \nonumber
u'' + u^2 = \alpha f(x)
\eea
for parameter $\alpha$ are studied. The problem is defined on the $x$ line with decay of both the solution $u(x)$ and the imposed forcing $f(x)$ as $|x|\to \infty$. The 
rate of decay of $f(x)$ is important and has a strong influence on the structure of the solution space. Three particular forcings are examined primarily: a rectilinear top-hat, a Gaussian, and a Lorentzian, the latter two exhibiting exponential and algebraic decay, respectively, for large $x$. The problem for the top hat can be solved exactly, but for the Gaussian and the Lorentzian it must be computed numerically in general. Calculations suggest that an infinite number of solution branches exist in each case. For the top-hat and the Gaussian the solution branches terminate
at a discrete set of $\alpha$ values starting from zero. A general asymptotic description of the solutions near to a termination point is constructed that also provides information on the existence of local fold behaviour. The solution branches for the Lorentzian forcing do not terminate in general. For large $\alpha$ the asymptotic analysis of Keeler, Binder \& Blyth (2018 `On the critical free-surface flow over localised topography',
{\it J. Fluid Mech.}, {\bf 832}, 73-96) is extended to describe the behaviour on any given solution branch using a method for glueing homoclinic connections.
\end{abstract}

%
%
%
\maketitle
%
%

\section{Introduction}

We investigate solutions to the nonlinear ordinary differential equation,
\begin{eqnarray}\label{maineq}
\frac{\dd^2 u}{\dd x^2} + u^2 = \alpha f(x),
\end{eqnarray}
for parameter $\alpha$, on the half-line $0 \leq x<\infty$ subject to the boundary conditions
\bea \label{mainbc}
\frac{\dd u}{\dd x}(0) = 0 \qquad \mbox{and} \qquad u(x) \to 0 \qquad \mbox{as} \qquad x\to \infty.
\eea
It is assumed that $f(0)=1$ and that $f\to 0$ as $|x|\to \infty$. The rate of decay for large $|x|$ is a delicate issue and has subtle and important implications for the solution. To highlight this feature of the problem, three particular forcing functions will be examined primarily: a top hat with compact support, a Gaussian and a Lorentzian, given by
\bea \label{forcings}
\hspace{-0.8in} f(x) = H(x+L)-H(x-L), \qquad f(x) = \mbox{e}^{-x^2} \qquad \mbox{and} \qquad f(x) = \frac{1}{1+x^2}
\eea
respectively, where $H(x)$ is the Heaviside function and $L$ is the half-width of the top hat.
Assuming that $f(x)=f(-x)$, as is the case for all of the forcings in \eqref{forcings}, the boundary condition \eqref{mainbc} may be viewed as providing an even solution over the entire $x$ line, and occasionally it will be helpful to discuss the problem in this context to illuminate some of the key features. 
Integrating \eqref{maineq} directly it is easily seen that
 \bea \label{solvability}
\alpha \int_0^\infty f(x)\,\dd x > 0
\eea
provides a necessary condition for a non-trivial solution to exist. For all three forcings in \eqref{forcings}
the integrand in \eqref{solvability} is non-negative and hence non-trivial solutions can only exist for $\alpha>0$. 

The problem is motivated by the study of free-surface flow of an inviscid, 
irrotational fluid over bottom topography. The forcing function $f(x)$ represents 
the negative of the topography so that the forcings in (\ref{forcings}) all correspond 
to a localised depression on an otherwise flat bottom. In the weakly-nonlinear limit of small 
forcing, the disturbance to the free-surface induced by the localised topography is 
described by the forced Korteweg-de Vries equation. The displacement of the free surface from its mean level is given by $u(x)$ governed by \eqref{maineq} assuming that the flow is steady and that the speed of the fluid far upstream of the depression is equal to the speed of small amplitude linear waves over a flat bottom (so that the Froude number for the flow 
is equal to unity).

\begin{figure}
\centering
\includegraphics[width=6.0in]{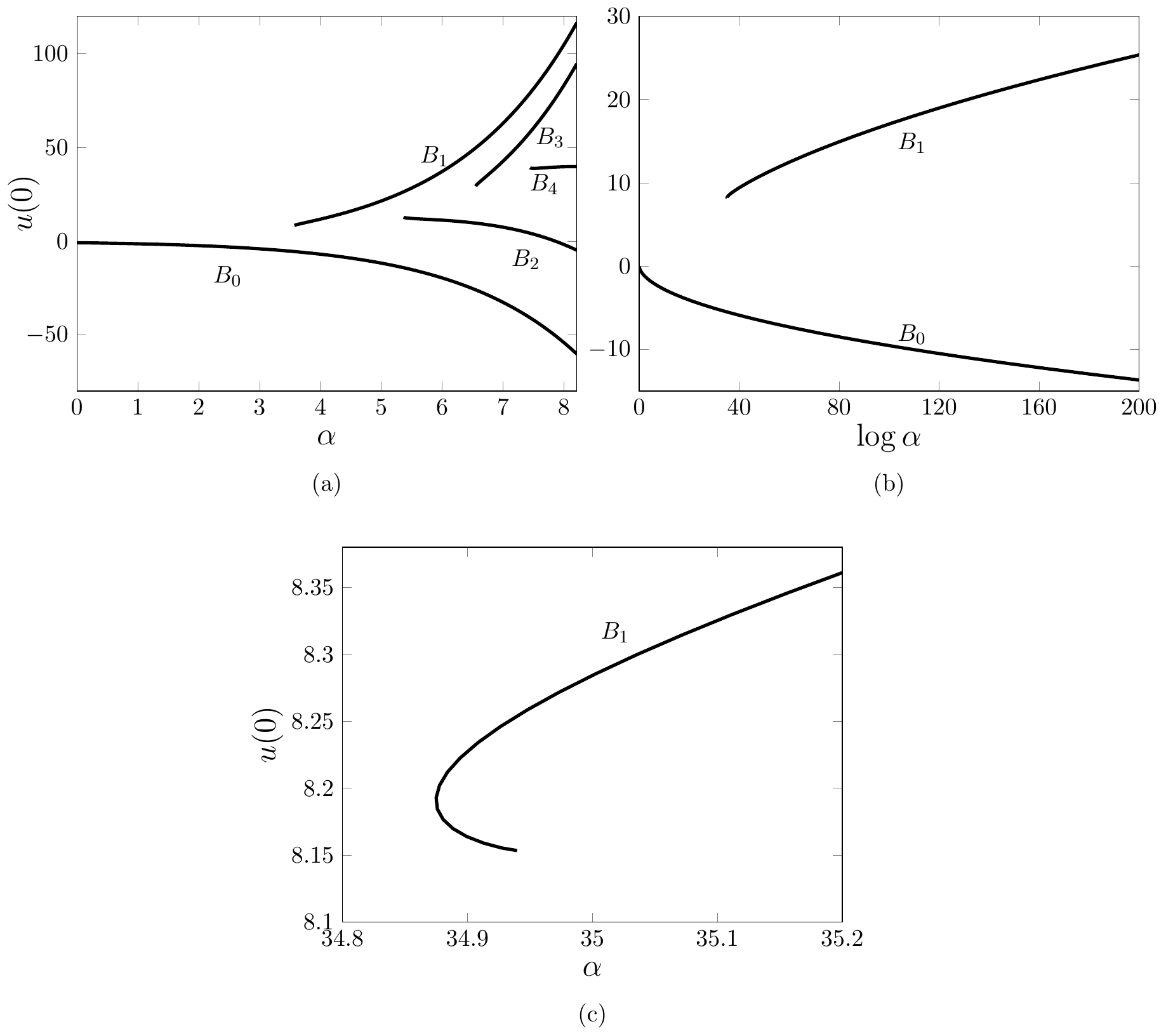}
\caption{Part of the solution space for the Gaussian forcing. The branches are labelled $B_n$ according to the integer number $n$ of local maxima in the solution $u(x)$.}
\label{fig:bifdig}
\end{figure}

Recently, Keeler {\it et al.} \cite{kbb} considered this problem for the Gaussian forcing. They made a number of observations about the solution space for $u(x)$ that require further mathematical explanation. In particular they 
presented numerical evidence that there exists an infinite number of distinct solution branches. To place the current work in context, figure \ref{fig:bifdig}\footnote{This is an adapted version of Figure 3(h,i) from Keeler, Binder, Blyth, On the critical free-surface flow over localised topography, {\it J. Fluid Mech.}, {\bf 832}, 73-96, used with permission.} shows part of the solution space uncovered by Keeler {\it et al.} \cite{kbb}, using $u(0)$ to characterise the solutions over a range of values of $\alpha$. In \cite{kbb} a traditional boundary layer analysis was used to construct asymptotic approximations both for small $\alpha$ and for large $\alpha$ that approximate the solutions on branch $B_0$ (in both limits) and on branch $B_1$ for large $\alpha$. The rest of the branches, labelled $B_n$ for integer $n$, are not captured by \cite{kbb}'s asymptotics and this provides one motivation for the present study. The present taxonomy for the solution branches differs from that used in \cite{kbb} and is motivated by the observation that the solution profiles on branch $B_n$ have $n$ local maxima. Since the solution spaces for all three of the forcings in (\ref{forcings}) share similar qualitative features (but with some key differences), the same taxonomy for the solution branches will be used in each case.
Keeler {\it et al.} \cite{kbb} provided solid but not conclusive numerical evidence that branch $B_1$ terminates at its leftmost end at a finite value of $\alpha$. The present work provides a deeper analysis of this issue.

The layout of the paper is as follows. In section \ref{sec:tophat} the top-hat forcing is considered; this problem has many of the important features also found for the smooth forcings but with the advantage that the solution can be found exactly. Next in section \ref{sec:farfield} the importance of the far-field decay rate for a smooth forcing is discussed, and the method for obtaining numerical solutions is described in section \ref{sec:numerics}. In section \ref{sec:gaussterm} an asymptotic analysis is presented that supports the termination of the branches $B_1$, $B_2$ etc. at finite $\alpha$ and an analysis that indicates that the branch $B_0$ terminates at $\alpha=0$. In section \ref{sec:agnesi} the case of a Lorentzian forcing is examined. Finally in section \ref{sec:glue} the method of homoclinic glueing is used to show how the large $\alpha$ solutions can be constructed for a general forcing with a local maximum. The appendices contain further details of the calculations for the homoclinic glueing, a Stokes line analysis for the Lorentzian forcing, and a discussion of a marginal case $f(x)=1/(1+x^4)$.

\section{Top hat forcing}\label{sec:tophat}

The top hat forcing, which takes the form given in \eqref{forcings}, provides an instructive model for the more technically challenging cases (the Gaussian and the 
Lorentzian forcings), not least because the solution can be obtained exactly in closed form.
\begin{figure}
	\centering
	\includegraphics[scale=0.7]{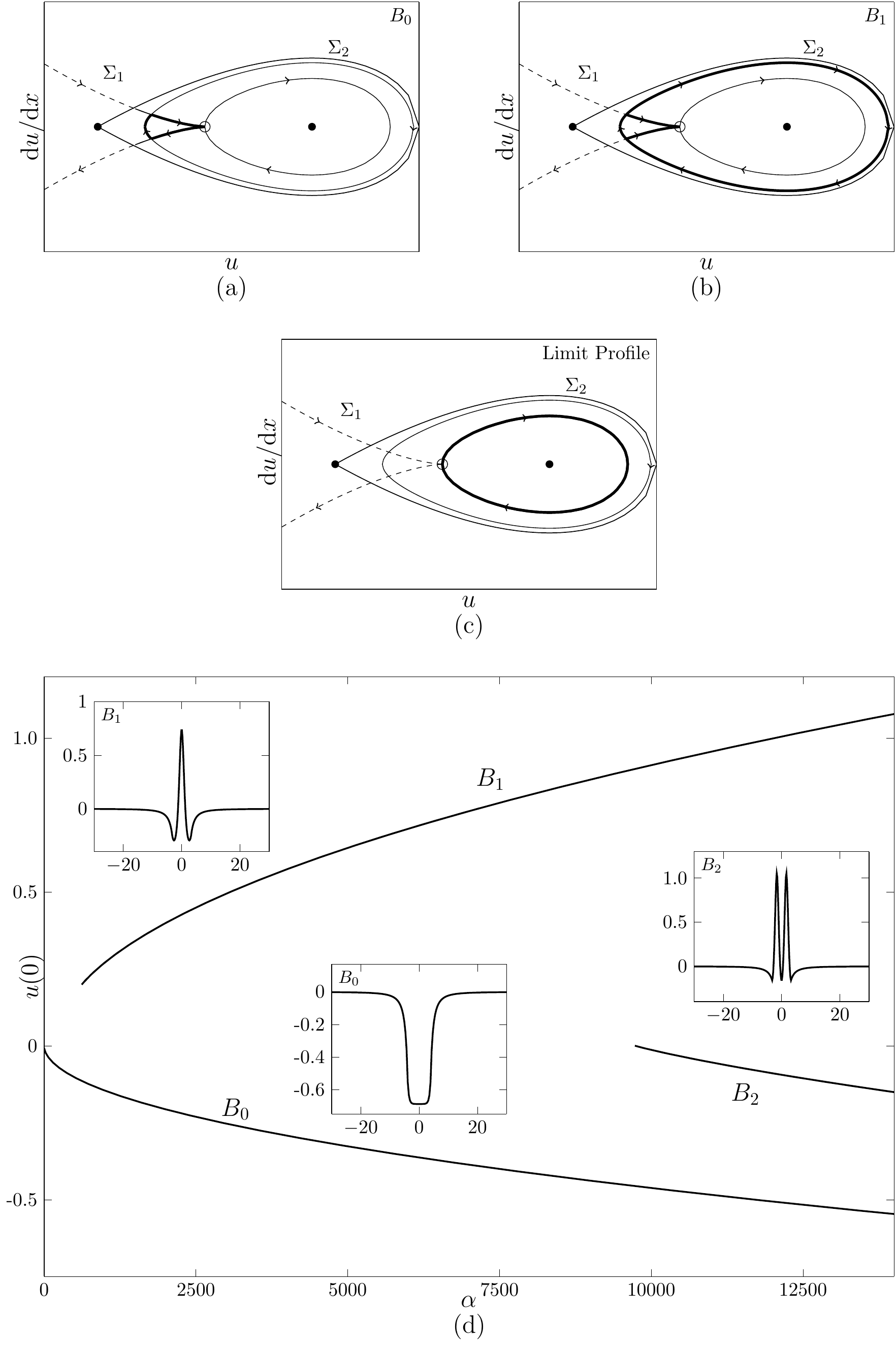}
	\caption{Top hat forcing: (a-c) phase portraits in the unforced ($\Sigma_1$) and forced ($\Sigma_2$) phase planes, and (d) the solution space. Panels (a-c) demonstrate the solution construction for branches $B_0$, $B_1$ and $B_2$ in the $(u,\dd u/\dd x)$ phase plane. The empty circles are located at the origin and correspond to a degenerate node in $\Sigma_1$, while the filled circles indicate the saddle point and centre in $\Sigma_2$. The broken lines in panels (a-c) are the stable and unstable manifolds for the degenerate node in $\Sigma_1$ and the thin solid lines are the orbits in $\Sigma_2$. A solution is indicated by a thick solid line. The solution branches are shown in panel (d) with insets showing sample solution profiles at the values $\alpha=20,618$ on $B_0$,  $\alpha=6,205$ on $B_1$, and $\alpha=13,584$ on $B_2$.}
	\label{fig:trench}
\end{figure}

A straightforward phase plane analysis nicely illustrates how the key features of the solution space emerge (see Binder \cite{binder} for a review of this technique applied to the KdV equation). 
The unforced phase plane, labelled $\Sigma_1$, corresponds to the homogeneous form of (\ref{maineq}) and is relevant outside of the top-hat's support where $|x|>L$.
It has a degenerate node at the origin, indicated in
figures~\ref{fig:trench}(a-c) by an empty circle, with a stable manifold and an unstable manifold on which
\bea
\tfrac{1}{2}u_x^2 + \tfrac{1}{3}u^3 = 0
\eea
holds and that are shown each with a broken line. 
The forced phase plane, labelled $\Sigma_2$, is relevant inside the top-hat support where $|x|<L$. It has a saddle point at $(u,\dd u/\dd x) = (-\alpha^{1/2},0)$ and a centre at $(\alpha^{1/2},0)$, both of which are shown in figures~\ref{fig:trench}(a-c) with filled circles. (Note that the phase portraits $\Sigma_1$ and $\Sigma_2$ are presented on the same scale.) Trajectories in $\Sigma_2$ satisfy
\bea
\tfrac{1}{2}u_x^2 + \tfrac{1}{3}u^3 = u + c
\eea
for constant $c$. These are shown with thin solid lines for different $c$ and comprise periodic orbits around the centre enclosed by a homoclinic orbit that connects the saddle to itself.
Solutions that satisfy the boundary conditions \eqref{mainbc} are indicated by thick solid lines in figures~\ref{fig:trench}(a-c).
In each case, starting from $x=-\infty$ the solution exits the origin and follows the unstable manifold in $\Sigma_1$
until $x=-L$ where it jumps instantaneously onto a periodic orbit in $\Sigma_2$. The trajectory jumps instantaneously back 
onto $\Sigma_1$ at $x=L$ and subsequently follows the stable manifold back into the origin as $x\to \infty$.
Thus the solutions are smooth everywhere except at $x=\pm L$ where the second derivative of $u$ is discontinuous.

Various possibilities arise while the trajectory is in $\Sigma_2$, depending on the value of $\alpha$. A solution that is 
negative-definite in $u$ can be constructed for any $\alpha>0$ by making only a partial excursion along the periodic orbit in the left-half plane of $\Sigma_2$, as is illustrated in figure~\ref{fig:trench}(a). Alternatively a trajectory may execute one cycle of the periodic orbit followed in general by a brief overshoot to make the connection back onto $\Sigma_1$, as is shown in figure~\ref{fig:trench}(b); however this is only possible if the top-hat is sufficiently wide and hence such a solution exists only when $\alpha>\alpha^*_1$, where $\alpha^*_1$ can be determined precisely and is given below. A countably infinite number of further options arises when $\alpha$ exceeds an increasing sequence of critical values, $\alpha=\alpha^*_n$, that can also be written down exactly. For each $n$ the solution executes $n$ cycles of a periodic orbit in $\Sigma_2$ followed by an overshoot to connect back to $\Sigma_1$. At the critical values themselves the solution executes exactly $n$ cycles along a periodic orbit in $\Sigma_2$, entering and leaving this plane from $\Sigma_1$ at the origin; in $\Sigma_1$ itself the solution is given by $u=0$ for all $x>|L|$. This critical case is illustrated in 
figure \ref{fig:trench}(c). 

%

The first three solution branches are shown in figure \ref{fig:trench}(d) together with some sample solution profiles.  
In all cases the phase plane trajectories are bounded within the homoclinic orbit in $\Sigma_2$; it follows that $-\alpha^{1/2}\leq u(x)<0$ on branch $B_0$ and $-\alpha^{1/2}\leq u(x)\leq 2\alpha^{1/2}$ on branches $B_n$ for $n=1,2,\cdots$. On the periodic orbit in $\Sigma_2$ for the critical case,
\bea
u(x) = (3\alpha)^{1/2}\,\mathrm{cn}^2\big((\alpha/3)^{1/4}x;1/2\big),
\eea
where $\mathrm{cn}$ is a Jacobi Elliptic function. Comparing the period of this form to the width of the top-hat we find
that
\bea
\alpha^*_n = 48 n^4 K^4(1/\sqrt{2}) \approx 567n^4,
\eea
for $n=1,2,\cdots$, where $K$ is the complete elliptic integral of the first kind.

\section{Far-field decay for smooth forcings}\label{sec:farfield}

The unforced, homogeneous form of (\ref{maineq}) has the general solution that decays at infinity,
\bea \label{algdec}
u_H(x) = -\frac{6}{(x+x_0)^2} 
\eea
for arbitrary constant $x_0$. Assuming that
\bea \label{fork}
f(x) = o(1/x^4) \quad \mbox{as} \quad x\to \infty,
\eea
the generic far-field behaviour of the solution is
\bea \label{cabin}
u\sim u_H(x) \qquad \mbox{as} \quad x\to \infty,
\eea
having a single degree of freedom, namely $x_0$ in \eqref{algdec}, which is effectively determined via the choice of $\alpha$.
The large $x$ balance between the first term on the left-hand side of (\ref{maineq}) and the forcing on the right-hand side,
\bea \label{astar}
\frac{\dd^2 u}{\dd x^2} \sim \alpha^*_n f(x) \qquad \mbox{as} \qquad x \to \infty
\eea
so that
\bea\label{astar2}
u\sim \alpha_n^*\int_x^{\infty} (x'-x)f(x')\,\dd x' \qquad \mbox{as} \qquad x \to \infty
\eea
is then also possible  but will occur only for certain special values of the parameter, 
$\alpha^*_n$, the behaviour \eqref{astar2} involving zero degrees of freedom (in integrating \eqref{astar} to obtain \eqref{astar2} both constants of integration must be fixed to ensure the far-field behaviour \eqref{mainbc} is satisfied).
Such a balance neglects the nonlinear term in (\ref{maineq}) and this is justified provided that \eqref{fork} holds.
If \eqref{fork} fails the generic far-field balance is between the nonlinear term and the forcing, given by
\bea \label{nadal}
u \sim - \alpha^{1/2}f^{1/2}(x),
\eea
where we have adopted the negative square root. The positive square root can be excluded on noting
that the linearisation at infinity $u = \pm \alpha^{1/2}f^{1/2}(x) + U(x)$ yields
\bea \label{chang}
\frac{\dd^2 U}{\dd x^2} \pm 2\alpha^{1/2}f^{1/2}U = 0.
\eea
If the positive square root is selected then a standard WKBJ analysis of \eqref{chang} yields the 
two linearly independent unbounded solutions
\bea
U \sim f^{-1/4}\,\mbox{exp}\big \{\pm \mathrm{i} \sqrt{2}\,\alpha^{1/4}\int^x f^{1/4}(x')\,\dd x'\big \}.
\eea
To exclude both of these requires two boundary conditions to be applied at infinity, but this leaves no freedom to enforce the symmetry condition at $x=0$ in \eqref{mainbc}. 
On the contrary, if the negative square root is selected, a single degree of freedom is retained since it is only necessary to exclude the exponentially growing solution to \eqref{chang}.
 The balance \eqref{nadal} occurs for the Lorentzian forcing, the third option in the list \eqref{forcings}, and this case will be examined in section \ref{sec:agnesi}.
 

\begin{figure}
\centering
\includegraphics[width=5.5in]{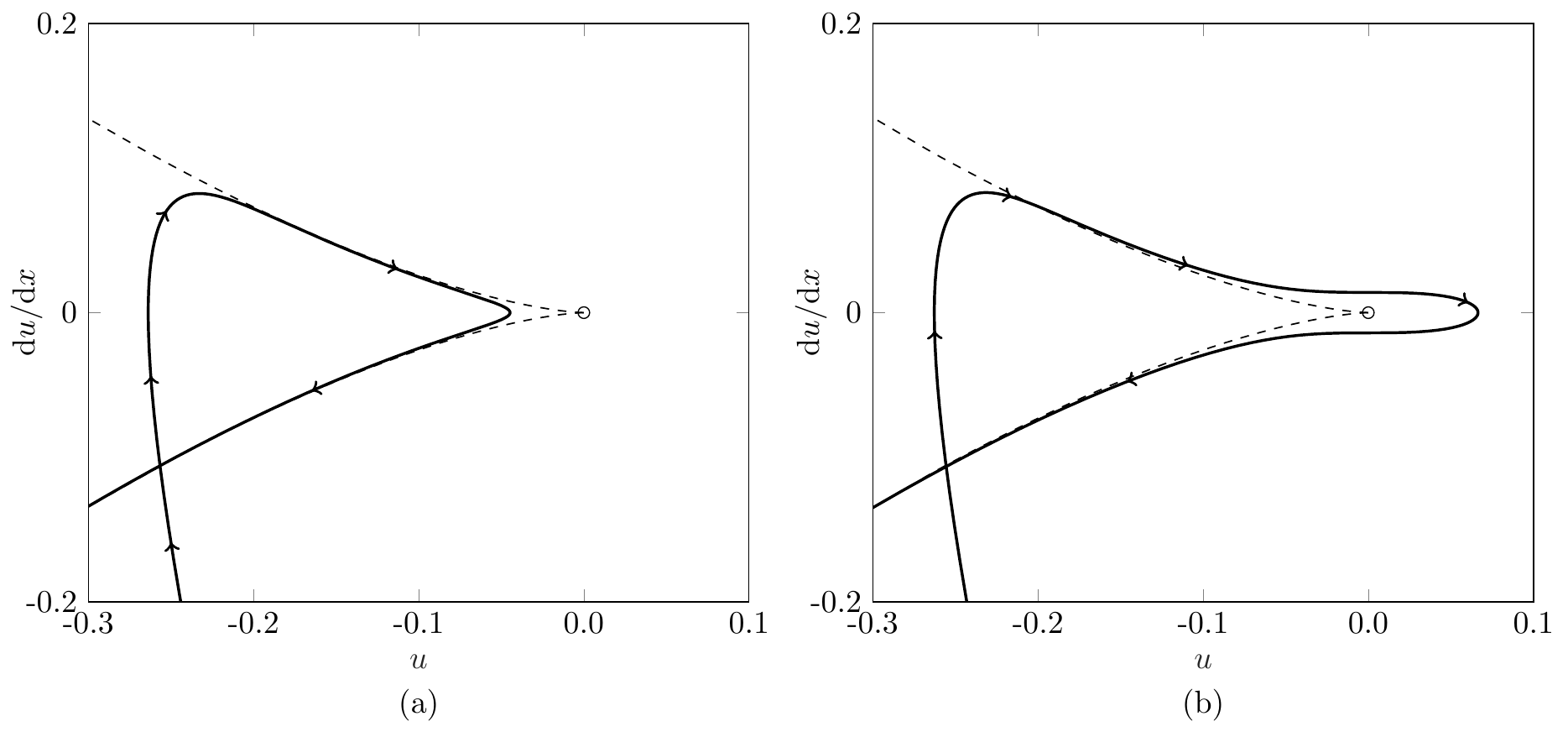}
\caption{Computed phase plane trajectories for the Gaussian forcing $f(x)=\mbox{exp}(-x^2)$ for $\alpha=36$ obtained by integrating \eqref{agfirst} with $u_0=8.5457$ for panel (a) and $u_0=8.5452$ for panel (b). Only a close-up near to the origin is shown. The empty circle and the broken lines correspond to the degenerate node and the stable/unstable manifolds in the unforced phase plane $\Sigma_1$.} 
\label{fig:shooting}
\end{figure}

\section{Numerical computation}\label{sec:numerics}

Numerical computations for the Gaussian forcing were carried out in \cite{kbb}. For any of the forcings in \eqref{forcings} it is expedient to first rewrite \eqref{maineq} as a first order system and then to solve the initial value problem
\bea \label{agfirst}
\bm{u}' = \bm{F}(\bm{u}), \qquad \bm{u}(0) = (u_0,0)^T,
\eea
where $\bm{u} = (u,\dd u/\dd x)^T$ and $\bm{F} = (\dd u/\dd x,\alpha f-u^2)^T$ for some $u_0$ to be found such that
$\bm{u}\to \bm{0}$ as $x\to \infty$ to fulfil \eqref{mainbc}. Thus a solution trajectory in the $(u,\dd u/\dd x)$ phase plane must 
ultimately enter the origin and, disregarding the degenerate behaviour \eqref{astar}, it must do so in the second quadrant.
For a Gaussian forcing, according to \eqref{cabin} it will enter the origin along the stable manifold
of the degenerate node in the unforced phase plane $\Sigma_1$ defined in section \ref{sec:tophat}. The computations for the Lorentzian forcing are particularly challenging as linearising about the far-field decay \eqref{nadal} by writing 
\bea \label{nadal2}
u(x) \sim - \frac{\alpha^{1/2}}{x} + U(x)
\eea
as $x\to \infty$ requires that
\bea
U'' - \frac{2}{x}U = 0
\eea
one solution of which,
\bea
U \propto x^{1/2}I_1\left(2\sqrt{2}x^{1/2}\right),
\eea
where $I_1$ is a modified Bessel function, grows exponentially for large $x$. 
Hence on shooting from $x=0$, any deviation from the required solution will rapidly grow.

\begin{figure}
\centering
\includegraphics[width=5.5in]{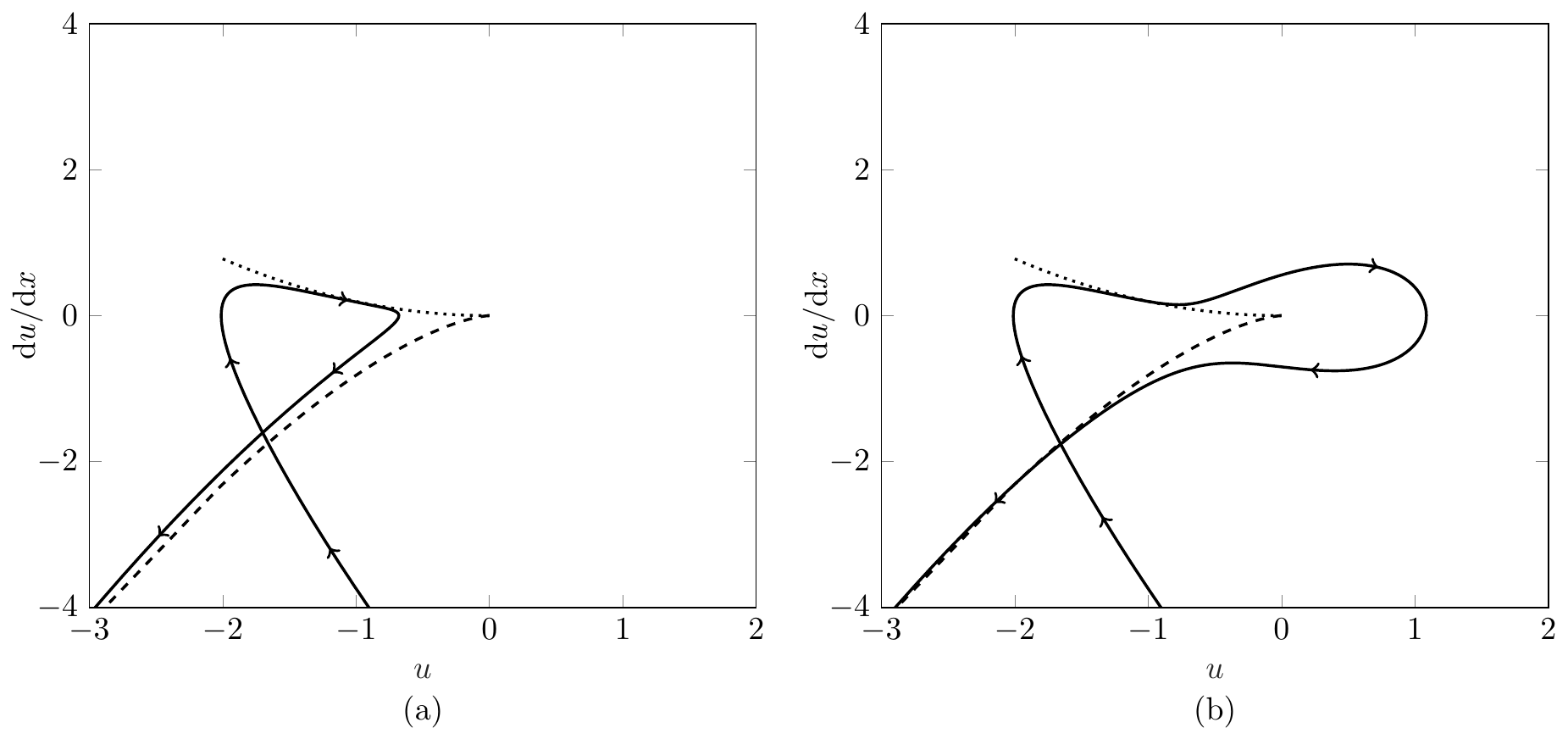}
\caption{Computed phase plane trajectories for the Lorentzian forcing $f(x)=1/(1+x^2)$ for $\alpha=26.44$ obtained by integrating \eqref{agfirst} with $u_0=8.298755$ for panel (a) and $u_0=8.298750$ for panel (b). Only a close-up near to the origin is shown. The empty circle and the thick broken line correspond to the degenerate node and the unstable manifold in the unforced phase plane $\Sigma_1$. The thin broken line indicates the behaviour $\dd u/\dd x = u^2/\alpha^{1/2}$ according to the leading order term in \eqref{nadal2}.} 
\label{fig:shooting_ag}
\end{figure}

In computational practice on a finite precision machine any choice for $u_0$ will result in `finite-time' blow up with $u\sim u_H$ as $x\to -x_0$ (with $x_0<0$), so that the phase plane trajectory converges to the unstable manifold in $\Sigma_1$. Nevertheless a solution can be detected to good accuracy by using a bisection approach. This is illustrated in figures \ref{fig:shooting} and \ref{fig:shooting_ag} for the Gaussian and Lorentzian forcings respectively. 
The trajectories were computed by integrating \eqref{agfirst} using the fourth-order Runge-Kutta method starting in each figure from two carefully selected positive values of $u_0$ (only a close-up near to the origin is shown). Since in both figures the two trajectories veer either side of the origin, assuming that $\bm{U}(x)$ depends continuously on $u_0$ there must exist a $u_0$ such that the corresponding trajectory reaches the origin and the far-field condition is satisfied.

Numerical calculations reveal that on a solution branch with a termination point the generic behaviour (\ref{algdec}) is found at all points along the branch except at the termination point where the singular far-field decay (\ref{astar}) is found. This is what was found, for example, in \cite{kbb} for a Gaussian forcing where at the termination point the decay is superexponential, corresponding to (\ref{astar}), and given by
\bea
u \sim \alpha_n^* \frac{\ee^{-x^2}}{4x^2} \qquad \mbox{as} \qquad x \to \infty.
\eea
Such solutions may be viewed as eigenmodes associated with eigenvalues $\alpha_n^*$, contrasting from solutions satisfying \eqref{cabin} in existing only for discrete values of $\alpha$ and exhibiting the maximal rate of decay as $x\to \infty$. We leave their relevance to applications as an open question.

%

\section{Branch termination}\label{sec:gaussterm}

The numerical calculations of \cite{kbb} for a Gaussian forcing suggest that branch $B_0$ terminates at $\alpha=0$ and branches $B_n$, $n=0,1,2,\ldots$ terminate at some value $\alpha_n>0$. In this section we present an asymptotic description of the branch termination in both cases.


\subsection{Termination at $\alpha=0$}\label{sec:smallalp}

As was noted in the Introduction non-trivial solutions to the problem \eqref{maineq}, \eqref{mainbc} for the stated class of forcing functions $f(x)$ exist only if $\alpha>0$. Therefore the solution branch $B_0$ that enters the origin in figure \ref{fig:bifdig} cannot pass into the left-half plane. Keeler {\it et al.} \cite{kbb} gave an asymptotic description of solutions on this branch for small $\alpha$. In fact branch $B_0$ must terminate at the origin. To demonstrate this, it is helpful to recapitulate some of the key details of the small $\alpha$ analysis.

\begin{figure}
  \centering
  \includegraphics[scale=0.7]{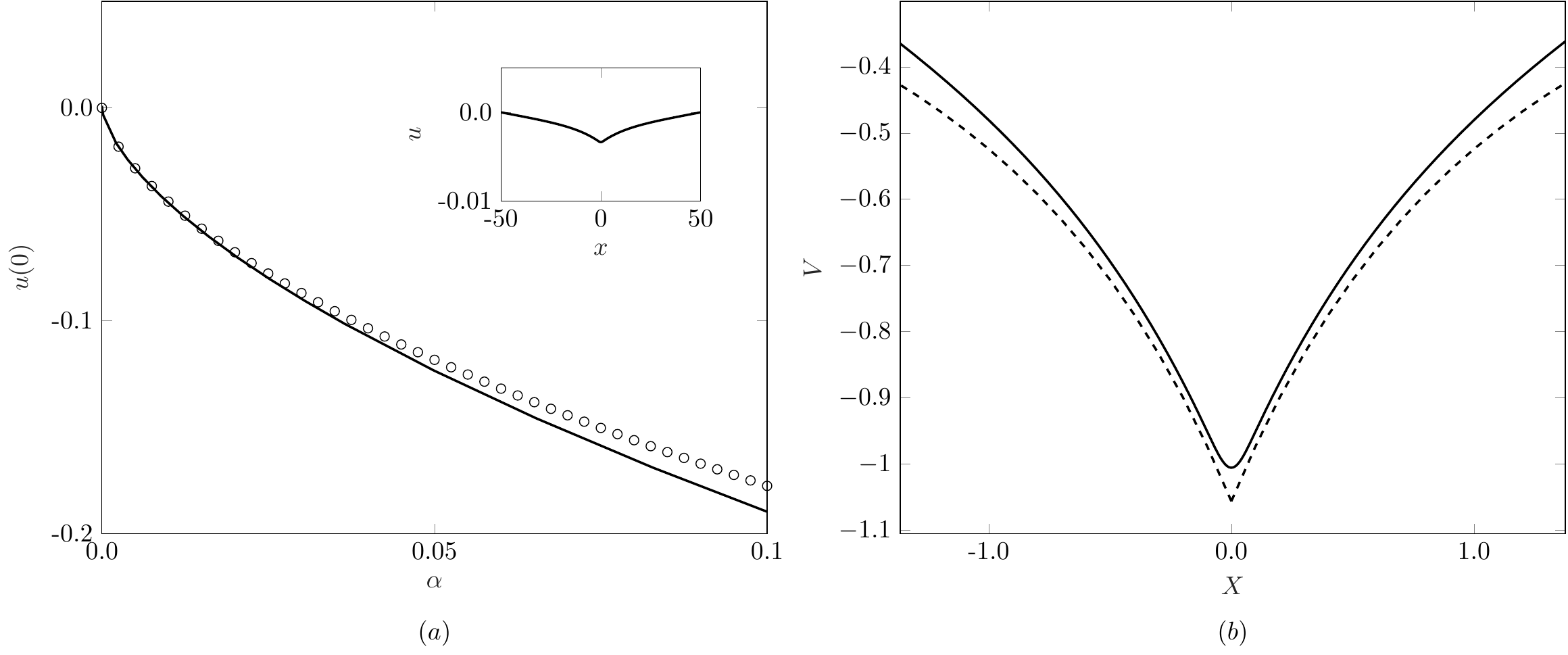}
  \caption{(a). The solution branch $B_0$ solution for the Gaussian forcing near to $\alpha=0$. The solid line is the numerically computed branch and the dotted markers represent the asymptotic approximation \eqref{ashes} with $M=\sqrt{\pi}$ with the next order correction included (see \cite{kbb}), namely $u(0)\sim -(3^{1/3}\pi^{2/3}/2)\alpha^{2/3} + \alpha/2$. The inset diagram is a solution profile on the outer scale when $\alpha = 1.928\times 10^{-4}$. (b). The same solution profile as for the inset in (a) but shown on the outer $(X,V)$ scale according to \eqref{saouter}. The solid line is the numerical solution and the dashed line is the asymptotic approximation \eqref{paine} with $X_0=2\cdot 3^{1/3}/\pi^{1/6}\approx 2.38$ given by \eqref{pooch}.}
  \label{fig:small_alpha}
\end{figure}

The expansion proceeds as $u(x;\alpha) = u(0;\alpha) + \alpha v(x)+\cdots$, noting that $u(0;\alpha)=O(\alpha^{2/3})$ which follows from the matching carried out below. Substituting into (\ref{maineq}), at leading order we obtain the linearised form
\begin{eqnarray}
\frac{\dd^2 v}{\dd x^2} = f(x),
\end{eqnarray}
with the boundary conditions $v(0) = \dd v/\dd x(0) = 0$. The far-field behaviour $v\sim Mx/2$ as $x\to \infty$, where the mass
\bea \label{mint}
M = 2\int_0^\infty f(x) \, \dd x,
\eea
suggests the outer scaling on which the nonlinear term is restored,
\bea \label{saouter}
u = \alpha^{2/3} V(X), \qquad x = \alpha^{-1/3}X. 
\eea
Under this scaling (\ref{maineq}) becomes
\bea \label{saV}
\frac{\dd^2 V}{\dd X^2} + V^2 = \alpha^{-1/3}f(\alpha^{-1/3}X)
\eea
Assuming $f(x) = o(1/x)$ as $x\to \infty$ the right hand side of (\ref{saV}) vanishes to leading order (note that this condition on $f$ is also required for the integral in (\ref{mint}) to be defined), and the solution subject to $V\to 0$ as $X\to \infty$ is 
\bea
V = -\frac{6}{(X_0 + |X|)^2}
\label{paine}
\eea
as shown in figure~\ref{fig:small_alpha}(b) for constant $X_0$. (Modulus bars have been included in \eqref{paine} to highlight the singular behaviour at $X=0$ discussed in detail below.) This behaves as 
\bea \label{gromit}
V \sim - \frac{6}{X_0^2} + \frac{12|X|}{X_0^3} + \cdots
\eea
as $X \to 0$.
Matching with the solution on the inner scale yields
\bea \label{pooch}
X_0 = 2(3/M)^{1/3}, 
\eea
and
\bea
u(0;\alpha) \sim -\left(\frac{3^{1/3}M^{2/3}}{2}\right)\,\alpha^{2/3}
\label{ashes}
\eea
as $\alpha \to 0^+$, which agrees with the leading order prediction in \cite{kbb} for a Gaussian forcing. This approximation is shown with the numerical solution in figure~\ref{fig:small_alpha}(a).

The boundary value problem written in the outer scalings (\ref{saouter}) best illustrates the way in which the branch terminates: 
the limit case $\alpha=0$ in (\ref{saV}) then corresponds to a na\"ive but natural replacement of 
(\ref{saV}) by 
\bea \label{saV2}
\frac{\dd^2 V}{\dd X^2} + V^2 = M\delta(X)
\eea
in the sense that (\ref{gromit}) and (\ref{pooch}) imply
\bea \label{massjump}
\left[\frac{\dd V}{\dd X}\right]^{0^+}_{0^-} = M.
\eea
(given the nonlinearity of (\ref{saV2}) this interpretation should of course be treated with considerable caution as will be seen in the next subsection). Thus the limit profile contains a corner 
and the branch cannot be continued.

\subsection{Termination at finite $\alpha$}\label{sec:term_finite}

As was discussed in section \ref{sec:farfield} the branches $B_n$ for $n=1,2,\ldots$ terminate at the special values $\alpha^*_n$. In this subsection a local asymptotic analysis is presented that describes the termination of an individual branch. To this end we write
\bea \label{pmsign}
\alpha = \alpha^*_n \pm \epsilon
\eea
with $0< \eps \ll 1$, where $\alpha^*_n$ is one of the special values discussed in section \ref{sec:farfield} at which the balance (\ref{astar}) holds and its value must be determined numerically. We therefore make the implict assumption that the forcing satisfies (\ref{fork}). The choice of sign in \eqref{pmsign} will be discussed below.

Introducing the expansion
\bea
u = u_0(x) + \eps u_1(x) + \cdots
\eea
and substituting into (\ref{maineq}) we obtain at leading order in $\eps$,
\bea
\frac{\dd^2 u_0}{\dd x^2} + u_0^2 = \alpha^*_n f(x)
\eea
with boundary conditions
\bea
\frac{\dd u_0}{\dd x}(0) = 0, \qquad u_0(x) \to 0 \quad \mbox{as} \quad x\to \infty.
\eea
By the definition of $\alpha^*_n$, the far-field decay of $u_0$ satisfies \eqref{astar}.
At first order we find
\bea \label{bvp1}
\frac{\dd^2 u_1}{\dd x^2} + 2u_0 u_1 = \pm f(x)
\eea
with 
\bea \label{bvp1bc}
\frac{\dd u_1}{\dd x}(0) = 0, \qquad u_1(x) = o(x) \quad \mbox{as} \quad x\to \infty.
\eea
The latter of these conditions is required for the matching to be described below and implies that
\bea
u_1(x) \to \pm \kappa \quad \mbox{as} \quad x\to \infty,
\eea
where the constant $\kappa$ is determined as part of the solution to the boundary value problem (\ref{bvp1}), (\ref{bvp1bc}).

An intermediate region holds where $f(x) = O(\epsilon)$ and $u(x) = O(\eps)$, and where the particular scaling on $x$ depends on the form of $f(x)$. However, the balance (\ref{astar}) still holds in this region and, fortunately, the linearity of the relation (\ref{astar}) implies that the solution can be expressed in the form $u=\eps v$, where $v$ is a linear combination of the constant $\kappa$ and the far-field form of $u_0$.
The latter becomes negligible outside of this region and on the outer scaling where $x=\eps^{-1/2}X$, with $X=O(1)$, writing $u=\eps v_0(X) + \cdots$ requires at leading order that
\bea \label{finite_outer}
\frac{\dd^2 v_0}{\dd X^2} + v_0^2 = 0.
\eea
The solution that decays in the far-field is
\bea \label{v0sol}
v_0 = -\frac{6}{\left( |X| + (6/|\kappa|)^{1/2}  \right)^2},
\eea
where we have included modulus bars to highlight the singular behaviour at $X=0$ to be discussed below.  Associated with the solution (\ref{v0sol}) is the requirement that $\kappa<0$ if the plus sign in (\ref{pmsign}) is used, meaning that near to the termination point the branch is such that $\alpha>\alpha^*_n$, and the requirement that $\kappa>0$ if the minus sign in (\ref{pmsign}) is used so that $\alpha<\alpha^*_n$ local to the termination point. 
These requirements ensure a match with the solution on the inner scale.
If follows that a sufficient condition for the existence of a fold in the solution branch is that $\kappa>0$. We may infer from the numerical results shown in Figure \ref{fig:bifdig} that $\kappa>0$ for the Gaussian forcing. 
It should be emphasised that while the present analysis gives a self-consistent description of the behaviour close to a termination point, it does not prove their existence even for forcings that satisfy the far-field decay condition (\ref{fork}). It may preclude the existence of a termination point, however, if (\ref{fork}) is not satisfied. Key to the latter remark is the existence of two possible large $x$ balances  for forcings that satisfy (\ref{fork}), namely \eqref{cabin} and \eqref{astar2}, on which the analysis presented in this subsection depends. Forcings that do not satisfy (\ref{fork}) have only one possible large $x$ balance as discussed in section \ref{sec:farfield}.

Since $u(0) \sim u_0(0) + \eps u_1(0) + \cdots$ the branches in figure \ref{fig:bifdig} are locally linear and enter the termination points with finite slope. In common with the small $\alpha$ case discussed in section \ref{sec:smallalp} the limit profile at each termination point, on the outer scale (\ref{v0sol}), has a corner at $X=0$ with the jump in slope
\bea \label{jumper}
\left[ \frac{\dd v_0}{\dd X} \right]_{0^-}^{0^+} = 2\sqrt{\frac{2}{3}} \, |\kappa|^{3/2}.
\eea
Notably, and in contrast to (\ref{massjump}), the jump is not given simply in terms of the mass of $f(x)$ but relies on the numerically determined constant $\kappa$. This underscores the danger, alluded to in the previous section, of a na\"ive replacement of the right hand side of the problem on the outer scale, here \eqref{finite_outer}, with $M\delta(X)$, where $M$ is the mass of the forcing given in \eqref{mint}, since the jump in slope at $X=0$  depends on the solution of the problem on the inner scale according to \eqref{jumper}.

The preceding remarks can be placed on a firmer footing by noting that 
on taking the limit $\alpha\to 0$ in \eqref{saV} the right hand side formally approaches a delta function (e.g. Stakgold and Holst \cite{stak}, Theorem 2.2.4).
On the contrary, writing \eqref{maineq} for the outer scaling results in \eqref{finite_outer} with
$$
\alpha \epsilon^{-2}f(\epsilon^{-1/2}X)
$$
on the right hand side, and this does not approach a delta function in the limit $\epsilon\to 0$. Consequently, 
branch $B_0$ can be described by na\"ively replacing the right hand side of \eqref{maineq} with $M\delta(x)$ and then
following the type of phase plane analysis reviewed by Binder \cite{binder}, but the remaining branches
$B_n$ for $n\geq 1$ cannot be described in this way.


\section{Lorentzian forcing}\label{sec:agnesi}

\begin{figure}
\centering
\includegraphics[width=5.0in]{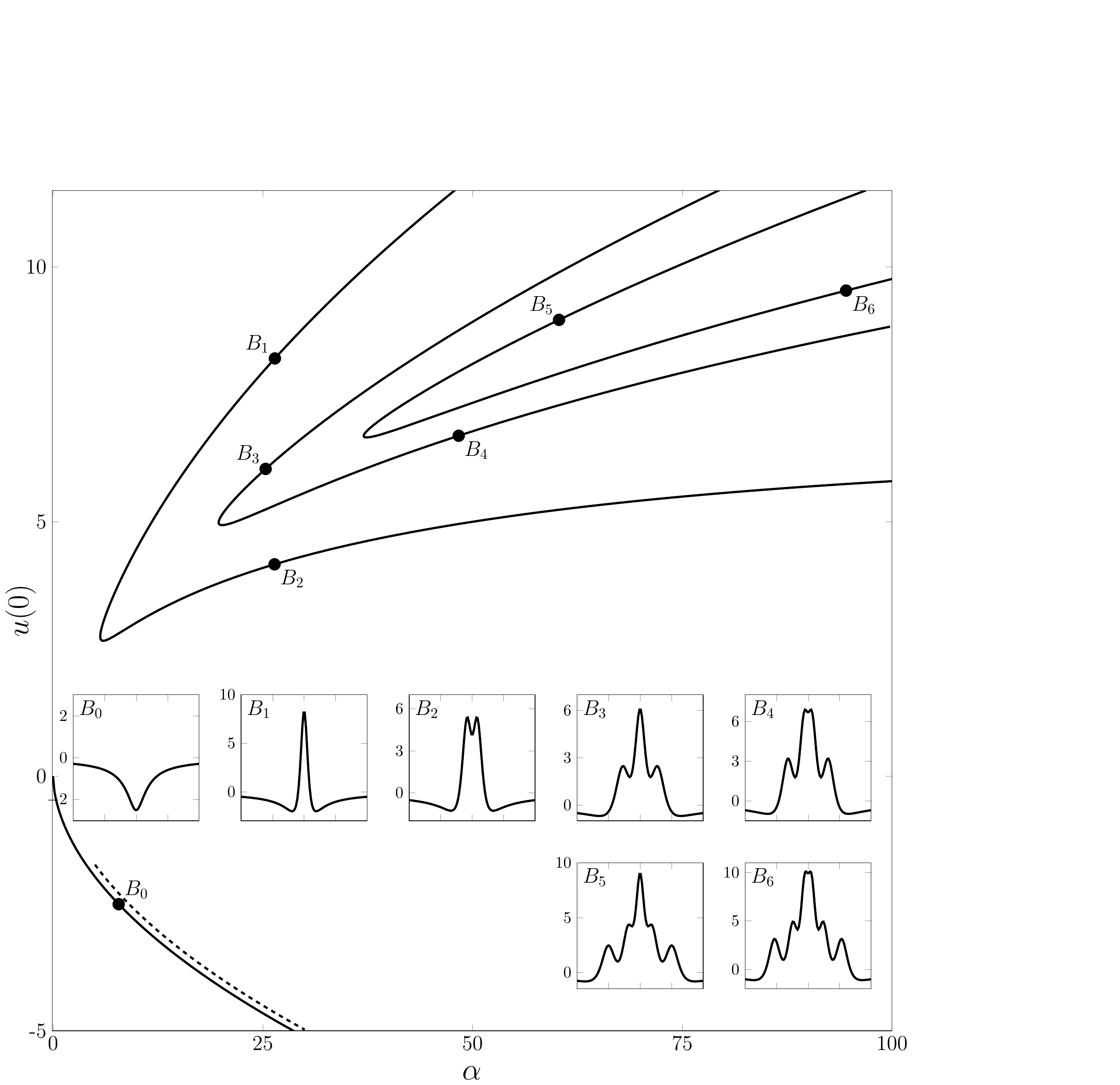}
\caption{Solution branches for the Lorentzian forcing, $f(x)=1/(1+x^2)$. All inset profiles are plotted in the range $x\in[-10,10]$ and correspond to the appropriate marker on the solution branches with the values $\alpha= 7.83$ ($B_0$), $\alpha= 26.44$ ($B_1$), $\alpha= 26.40$ ($B_2$), $\alpha= 25.34$ ($B_3$), $\alpha= 48.35$ ($B_4$), $\alpha= 60.31$ ($B_5$), $\alpha = 94.52$ ($B_6$).}
\label{fig:BVP_agnesi1}
\end{figure}

The numerically computed solution space for the Lorentzian forcing in \eqref{forcings} is similar in structure to that found for the Gaussian (see figure \ref{fig:bifdig}) with the crucial difference that the higher order branches $B_1$, $B_2$, etc. do not terminate at finite $\alpha$. As for the Gaussian forcing there is a $B_0$ branch of negative-definite solutions that exists for all positive $\alpha$ and which terminates at $\alpha=0$ as described in section \ref{sec:smallalp}. We do not expect to find branches that terminate at non-zero $\alpha$ for reasons discussed in section \ref{sec:farfield}. In fact we find that the branches $B_1$ and $B_2$ are in this case connected continuously as can be seen in Figure \ref{fig:BVP_agnesi1}. The insets in this figure show typical solution profiles some way along the upper and lower parts of the branch. Note that although the comment at the end of section \ref{sec:farfield} leaves open the possibility that eigenmode solutions exist for special values of $\alpha$ (as in the Gaussian case) corresponding to the choice of the plus sign in \eqref{chang}, our numerical computations suggest that such solutions do not, in fact, exist.

Following the success of the boundary-layer analysis of \cite{kbb} on branch $B_0$ for the Gaussian forcing, we are motivated to attempt a similar large $\alpha$ analysis for the Lorentzian, and this is considered in the following subsection. As for the Gaussian such an analysis cannot capture the higher order branches. These will be discussed in \ref{sec:glue}.

\subsection{Asymptotic approximation when $\alpha\gg 1$}\label{sec:agnesi_large}

For large $\alpha$ we rescale by writing $u=\alpha^{1/2}W(x)$ so that (\ref{maineq}) becomes
\bea
\asm \frac{\dd^2 W}{\dd x^2} + W^2 = \frac{1}{1+x^2},
\label{smith}
\eea
where $\mu=\alpha^{-1/2}$. We seek an asymptotic expansion in the form
\begin{equation}
W(x) = W_0(x) + \asm W_1(x) + \asm^2 W_2(x) +\cdots
\label{gayle}
\end{equation}
Substituting into (\ref{smith}) and working at successive orders, we obtain
\bea \label{cricket}
\hspace{-0.65in} W_0 = \pm \frac{1}{(1+x^2)^{1/2}}, \quad W_1 = -\frac{1}{2}\frac{(2x^2-1)}{(1+x^2)^2}, \quad
W_2 = \mp \frac{5}{8}\frac{(4x^4-20x^2+3)}{(1+x^2)^{7/2}}.
\eea
In general for $n\geq 1$,
\bea
W_n = -(2W_0)^{-1}\left( \sum_{k=1}^{n-1}W_kW_{n-k} + \frac{\dd^2 W_{n-1}}{\dd x^2}\right),
\label{rashid}
\eea
and it is straightforward to show that $W_n(x) = p_{2n}(x)/(1+x^2)^{(1+3n)/2}$, where $p_{2n}(x)$ is a polynomial of degree $2n$. For any $n$ an important observation is that $W_n$ satisfies the boundary conditions (\ref{mainbc}) irrespective of the choice of sign in (\ref{cricket}). This is reminiscent of the problems discussed by Chapman {\it et al.} \cite{CKA} and references therein, whereby for some ordinary differential equations with a small parameter, a simple asymptotic solution can be constructed that satisfies the required boundary conditions at any order and for all $x$ when no such solution to the problem in fact exists.

The difficulty can be traced to the fact that the expansion disorders in the neighbourhood of
the singularities in the leading order term in the 
asymptotic expansion ($W_0$ here) extended into the complex plane. If a Stokes line emanating from one of these 
singularities crosses the real axis, an exponentially small `beyond all orders' term is in general switched on at the point of crossing 
and this term eventually grows to corrupt the original expansion. In Appendix B we provide the relevant Stokes line 
analysis for the present problem. The conclusion is that if the plus sign in \eqref{cricket} is chosen then the expansion 
\eqref{gayle} is corrupted in the manner described. For the minus sign no such difficulty arises and the asymptotic 
expansion (\ref{gayle}) is valid for all $x$. 

Selecting the minus sign in (\ref{gayle}), we find $u(0) \sim -\alpha^{1/2} + 1/2$. This approximation is shown with a broken line in Figure \ref{fig:BVP_agnesi1}. As this figure suggests, further possibilities arise for the large $\alpha$ asymptotics in which \eqref{gayle} is regarded as an outer expansion to be matched to an inner boundary-layer solution around $x=0$ that describes a cluster of localised waves. These localised waves may be viewed as the connection of individual solitary-wave structures that are each described by a homoclinic orbit in an appropriately defined phase space, as will be discussed in the following section.  

\section{Homoclinic glueing}\label{sec:glue}

In this section we aim to describe asymptotic forms that approximate the solutions in the limit of large $\alpha$ for a fairly general class of forcing functions $f$.
It will be convenient to think in terms of symmetric solutions to \eqref{maineq} that are defined on the whole of the real line and that decay as $|x|\to \pm \infty$. To motivate the construction we rewrite (\ref{maineq}) on a boundary-layer scale 
by introducing the new variable $x=\alpha^{-1/4} y$ and setting $u(x)=\alpha^{1/2}U(y)$ to 
obtain
\begin{eqnarray}\label{barnacles}
\frac{\dd^2 U}{\dd y^2} + U^2 = f(\dg y),
\end{eqnarray}
where $\dg=\alpha^{-1/4}\ll 1$. Expanding the solution as $U = U_0(y) + \dg^2 U_1(y)+\cdots$,
and replacing the right hand side by its Taylor expansion $f(\delta y) = f(0) + (1/2)\delta^2 f''(0) y^2 + O(\delta^3)$,
at leading order we find (since $f(0)=1$)
\bea \label{inkling}
\frac{\dd^2 U_0}{\dd y^2} + U_0^2 = 1.
\eea
This has three bounded solutions of interest,
\bea \label{tweak}
\hspace{-0.35in}
(\mathrm{i})\,\,\, U_0=1, \qquad (\mathrm{ii})\,\,\, U_0 =-1, \qquad (\mathrm{iii}) \,\,\, U_0 = U_0^H \equiv 3\mbox{sech}^2(y/\sqrt{2}) - 1,
\eea
that correspond to equilibria in the $(U_0,\dd U_0/\dd y)$ phase plane (i and ii) and
a homoclinic orbit connecting $(-1,0)$ to itself in the same plane (iii). It is symmetric about $y=0$ and has the property
\bea \label{peso}
U_0^H \sim  -1  + 12 \ee^{\mp \sqrt{2}y} \quad \mbox{as $y\to \infty$}.
\eea
Its graph in physical space has a classical solitary-wave shape (see, for example, Billingham and King \cite{bilk}), and this suggests representing the wave-like parts of the $B_n$  
branch solutions by a collection of these homoclinics. 

The strategy is as follows: for $n$ odd, we seek to glue together $n$ of the homoclinics (iii) in \eqref{tweak} via asymptotic matching; for $n$ even, $n-1$ homoclinics are glued together either side of a central region in which solution (ii) in \eqref{tweak} predominates. We consider odd and even numbers of homoclinics separately in the following subsections.
In both cases the analysis is predicated on the assumption that $f''(0)<0$, so that the forcing has a local maximum at the origin (or, by a suitable shift, at any $x$ location). 
This condition is fulfilled both by the Gaussian and the Lorentzian forcings.

\subsection{Odd number of homoclinics}\label{sec:oddhomoclinic}

Continuing with the analysis, at next order we have
\bea \label{kwazii}
\frac{\dd^2 U_1}{\dd y^2} + 2U_0^HU_1 = \tfrac{1}{2}f''(0)y^2.
\eea
Henceforth in this subsection it is assumed that $U_0$ is given by the homoclinic (iii) in \eqref{tweak}. 
The general solution to (\ref{kwazii}) is given in Appendix A in terms of a symmetric and an antisymmetric complementary 
function $\Phi_s$ and $\Phi_a$, and particular integral form $\Phi_p^{(2)}$ that satisfies $\Phi_p^{(2)}(0)=\dd \Phi_p^{(2)}/\dd y (0)=0$. 
The solution that satisfies the boundary conditions
\bea \label{kwazii_cond}
U_1(0) = U_{10}, \qquad \frac{\dd U_1}{\dd y}(0) = 0
\eea
is given by
\bea \label{shellington}
U_1(y) = U_{10}\,\Phi_s(y/\sq) + \tfrac{1}{2}f''(0)\Phi_p^{(2)}(y/\sq),
\eea
where $U_{10}$ is a constant that will be determined later. We note from (\ref{app:asym}) that
\bea
\Phi_s(y/\sq) \sim -\tfrac{1}{16}\ee^{\sq y}, \qquad 
\Phi_p^{(2)}(y/\sq) \sim \tfrac{1}{4}(\log 2)\, \ee^{\sq y}
\eea
as $y\to \infty$. 

For a single homoclinic we require 
\bea \label{tunip}
U_{10} = 2(\log 2)f''(0) 
\eea
to exclude the exponential growth as $y\to +\infty$. This is the case considered in \cite{kbb} for the Gaussian forcing. The solution for $U(y)$ is then matched to the solution on 
the outer scale, where $x=O(1)$\footnote{The problem on the outer scale was discussed for the Lorentzian forcing in section \ref{sec:agnesi}. For the Gaussian see \cite{kbb}.}. 
More generally, using the asymptotic forms given in Appendix A we have as the required matching condition
\bea \label{gupA}
\hspace{-0.5in}
U \sim -1 + 12 \ee^{-\sq y} + \dg^2 \left( -\tfrac{1}{4}f''(0)(y^2+1) + \Lambda_1\ee^{\sq y} +\cdots \right) + \cdots 
\eea
as $y\to +\infty$, where 
\bea \label{ray}
\Lambda_1 = \tfrac{1}{8}(\log 2)f''(0) -\tfrac{1}{16} U_{10}.
\eea 
According to \eqref{tunip} the single homoclinic has $\Lambda_1=0$.  If $\Lambda_1>0$ then a second homoclinic is initiated at
\bea \label{manta}
y = Y_1(\dg) + y_1, \qquad Y_1 = \frac{1}{\sq}\log \left( \frac{12}{\Lambda_1\dg^2} \right),
\eea
 whereupon \eqref{gupA} becomes
\bea \label{gupB}
\hspace{-0.65in}
U \sim -1 + 12 \ee^{\sq y_1} + \dg^2 \left( -\tfrac{1}{4}f''(0)\big [ (Y_1+y_1)^2 + 1 \big] + \Lambda_1\ee^{-\sq y_1} + \cdots \right) + \cdots 
\eea
as $y_1 \to - \infty$. A third homoclinic is initiated in $y<0$ to maintain symmetry.

\vspace{0.125in}
\noindent
{\sc Remark 1}: The shift $Y_1(\delta)$ has been chosen to make the correction term $\delta^2 \Lambda_1\ee^{\sq y}$ in \eqref{gupA} of $O(1)$ and so that the two exponential terms in \eqref{gupA} are effectively interchanged in \eqref{gupB} to effect the matching between the homoclinics. 

\vspace{0.125in}
\noindent
For $y_1=O(1)$ we have the expansion
\bea
U=U_0^H(y_1)+\dg^2 \theta_1(y_1; \log (1/\dg))+\cdots,
\eea
where $U_0^H$ was given in (\ref{tweak}) and 
\bea \label{dashi}
\frac{\dd^2 \theta_1}{\dd y_1^2} + 2U_0^H\theta_1 = \tfrac{1}{2}f''(0)(Y_1 + y_1)^2,
\eea
provided that $\delta^2 Y_1\ll 1$, a restriction that will be discussed in more detail below. 
Henceforth the subscripts on $\theta$ will label the homoclinic sequence rather than the $\dg$ expansion. Also, as is 
suggested by the notation, here and subsequently we shall lump all additional logarithmic factors into $\theta_i$ in order to 
obtain algebraic rather than simply logarithmic accuracy. 

\vspace{0.125in}
\noindent
{\sc Remark 2}: The inner expansion at $O(\dg^4)$ will lead to a $\dg^2 \ee^{\sq y_1}$ term 
in (\ref{gupB}). This will simply trigger a $\dd U_0^H/\dd y_1$ complementary function in the 
solution to (\ref{dashi}), which corresponds to an $O(\dg^2)$ translation in $y_1$. This
can be safely ignored since we do not seek to determine $O(\dg^2)$ corrections to the 
locations of the maxima.

\vspace{0.125in}
Inspecting (\ref{gupB}) and (\ref{dashi}) we decompose $\theta_1$ as
\bea \label{gupC}
\theta_1 = \Lambda_1\phi(y_1) + \tfrac{1}{2}f''(0)\big( Y_1^2 \psi_0(y_1) + 2 Y_1\psi_1(y_1) + \psi_2(y_1)\big),
\eea
where $\phi$ and the $\psi_k$, $k=0,1,2$ satisfy the problems stated in Appendix A. To perform the matching we demand that
\bea \label{octopod}
\phi \sim \ee^{-\sq y_1}
\eea
as $y_1 \to -\infty$, where the expansion does not include a term of the form $a\ee^{\sq y_1}$ for any constant $a$.
We also require that 
\begin{align}\label{cbeebies}
\psi_k \sim -\tfrac{1}{2}y_1^k  \quad (\mbox{for } k=0,1), \qquad \psi_2 \sim -\tfrac{1}{2}(y_1^2+1),
\end{align}
as $y_1\to -\infty$, where in both cases \eqref{cbeebies} the expansions do not include either $a\ee^{\sq y_1}$ or $a\ee^{-\sq y_1}$ for any constant $a$. The given stipulations for the $\cdots$ in both \eqref{octopod} and \eqref{cbeebies} are made to remove the translational invariance alluded to above to ensure a unique solution.

The solutions for $\phi$ and the $\psi_k$ are given in \ref{app:glue}; here it is sufficient to note that
\bea \label{gupE}
\theta_1 \sim \Lambda_2 \ee^{\sq y_1}
\eea
as $y_1\to \infty$, where $\Lambda_2 =  \Lambda_1 + (\sq/8)f''(0)Y_1$. Therefore a triple homoclinic solution requires that $\Lambda_2=0$, that is
\bea \label{manta2}
\hspace{-0.8in}
\Lambda_1 = -\tfrac{\sq}{8}f''(0)Y_1(\delta)
=
-\tfrac{1}{8}f''(0)\log \left( \frac{12}{\Lambda_1\dg^2} \right)
\sim -\tfrac{1}{4}f''(0)\log (1/\delta)
\eea
using \eqref{manta}. It follows from \eqref{ray} that $U_{10} \sim 4f''(0)\log (1/\delta)$.

If \eqref{gupE} is not satisfied then \eqref{gupA} is replaced as the matching condition to the next region by   
\begin{align} \label{gupF}
\nonumber
U \sim -1 + 12 \ee^{-\sq y_1} + \dg^2 \Big( & -\tfrac{1}{4}f''(0)(y_1^2+1 + 2Y_1y_1 + Y_1^2)  \\
& + \Lambda_2 \ee^{\sq y_1} +\cdots \Big ) + \cdots 
\end{align}
as $y_1\to +\infty$. Hence for $\Lambda_2>0$ a fourth homoclinic arises for $y_1=Y_2-Y_1+y_2$ where 
\bea
Y_2(\delta) = \frac{1}{\sq}\log \left(\frac{12}{\Lambda_2\delta^2}\right) + Y_1(\delta),
\eea
a fifth homoclinic being initiated in $y<0$ by symmetry, and
\eqref{gupF} implies the matching condition
\begin{equation} \label{gupG}
\nonumber
U \sim -1 + 12 \ee^{\sq y_2} + \dg^2 \Big[  -\tfrac{1}{4}f''(0)\Big((Y_2+y_2)^2+1\Big)  
 + \Lambda_2 \ee^{-\sq y_2} +\cdots \Big ] + \cdots 
\end{equation}
as $y_2\to -\infty$. In the new homoclinic region we therefore write
\bea
U = U_0^H(y_2)+\dg^2 \theta_2(y_2; \log (1/\dg))+\cdots,
\eea
where 
\bea \label{gupH}
\theta_2 = \Lambda_2\phi(y_2) + \tfrac{1}{2}f''(0)\big( Y_2^2 \psi_0(y_2) + 2 Y_2\psi_1(y_2) + \psi_2(y_2)\big),
\eea
similar to \eqref{gupC}. Thus the sequence is now established with
\bea
\Lambda_{n+1} = \Lambda_n + \frac{\sq}{8}f''(0)Y_n, \qquad Y_{n+1} = \frac{1}{\sq}\log\left(\frac{12}{\Lambda_{n+1}\delta^2}\right) + Y_n,
\eea 
and $y_n = Y_{n+1}-Y_n + y_{n+1}$ for $n=1,2,\cdots$. 
The homoclinic sequence  for each successive $n$ corresponds to the value $U_{10}$, given in \eqref{ray}, such 
that $\Lambda_n=0$. The sequence has a homoclinic at $y=0$ and, when $n\geq 2$, at $y=\pm Y_k$ for $k=1,\cdots n-1$. At leading order in $\log (1/\delta)$,
\bea \label{bluh}
Y_n\sim \sq \,n \log (1/\delta).
\eea

\subsection{Even number of homoclinics}

The first pair of homoclinics is located at $y=\pm \ye_1(\delta)$, where $\ye_1$ is to be found (note that $Y_1$ now differs from that given in section \ref{sec:oddhomoclinic}). Sufficiently close to $y=0$ the 
expansion
\bea \label{gupI}
U = -1 + \delta^2 \Ue_1(y) + \cdots
\eea
holds, where the leading order term corresponds to (ii) in \eqref{tweak}. Note that the choice (i) in \eqref{tweak} was ruled out in \cite{kbb} and is ruled out here for the same reason. At first order
\bea
\frac{\dd \Ue_1}{\dd y} - 2 \Ue_1 = \tfrac{1}{2}f''(0)y^2
\eea
with solution
\bea \label{gupJ}
\Ue_1 = -\tfrac{1}{4}f''(0)(y^2+1) + \frac{12}{\delta^2}\left(\ee^{\sq(y-\ye_1)} + \ee^{-\sq(y+\ye_1)}\right),
\eea
where the constants of integration have been set to ensure a match with the homoclinics at $y=\pm \ye_1$. 

\vspace{0.125in}
\noindent
{\sc Remark 3}: \eqref{gupJ} appears to be inconsistent with the expansion \eqref{gupI}; in fact $\ye_1$ will be such that 
the $1/\delta^2$ terms in \eqref{gupJ} are of $O(1/\delta)$ approximately, a statement that will be made more precise below, so $|\delta^2 U_1| = O(\delta)\ll 1$.

\vspace{0.125in}
\noindent
Setting $y=\ye_1+\hy_1$, and inserting \eqref{gupJ}, \eqref{gupI} becomes
\bea
\hspace{-0.8in}
U = -1 + 12\ee^{\sq \hy_1} + \delta^2\left(  -\tfrac{1}{4}f''(0)\big( (\ye_1+\hy_1)^2+1\big) + \frac{12}{\delta^2}\ee^{-2\sq\ye_1}\ee^{-\sq \hy_1}\right) + \cdots,
\eea
which motivates writing $u = U_0^H(\hy_1) + \delta^2 \theta_1(y_1;\log(1/\delta))$ for $\hy_1=O(1)$, where
\bea \label{humma}
\theta_1 =  \frac{12}{\delta^2}\ee^{-2\sq \ye_1}\phi(\hy_1) + \tfrac{1}{2}f''(0)\left(\ye_1^2 \psi_0(\hy_1) + 2\ye_1\psi_1(\hy_1) + \psi_2(\hy_1)\right)
\eea
(the functions $\phi$ and $\psi_k$, for $k=0,1,2$, were defined in section \ref{sec:oddhomoclinic}). Once $\ye_1$ is determined below, it can be confirmed {\it a posteriori} that each of the terms in \eqref{humma} are at worst logarithmic in $\delta$. 
Using the details given in \ref{app:glue},
\bea \label{gupK}
\theta_1 \sim \Lambda_1\ee^{\sq \hy_1}, \qquad \quad \Lambda_1 = \frac{12}{\delta^2}\ee^{-2\sq \ye_1} + \frac{\sq}{8}f''(0)\ye_1,
\eea
as
$\hy_1\to +\infty$. 

%
\begin{figure}
\centering
\includegraphics[width=5.5in]{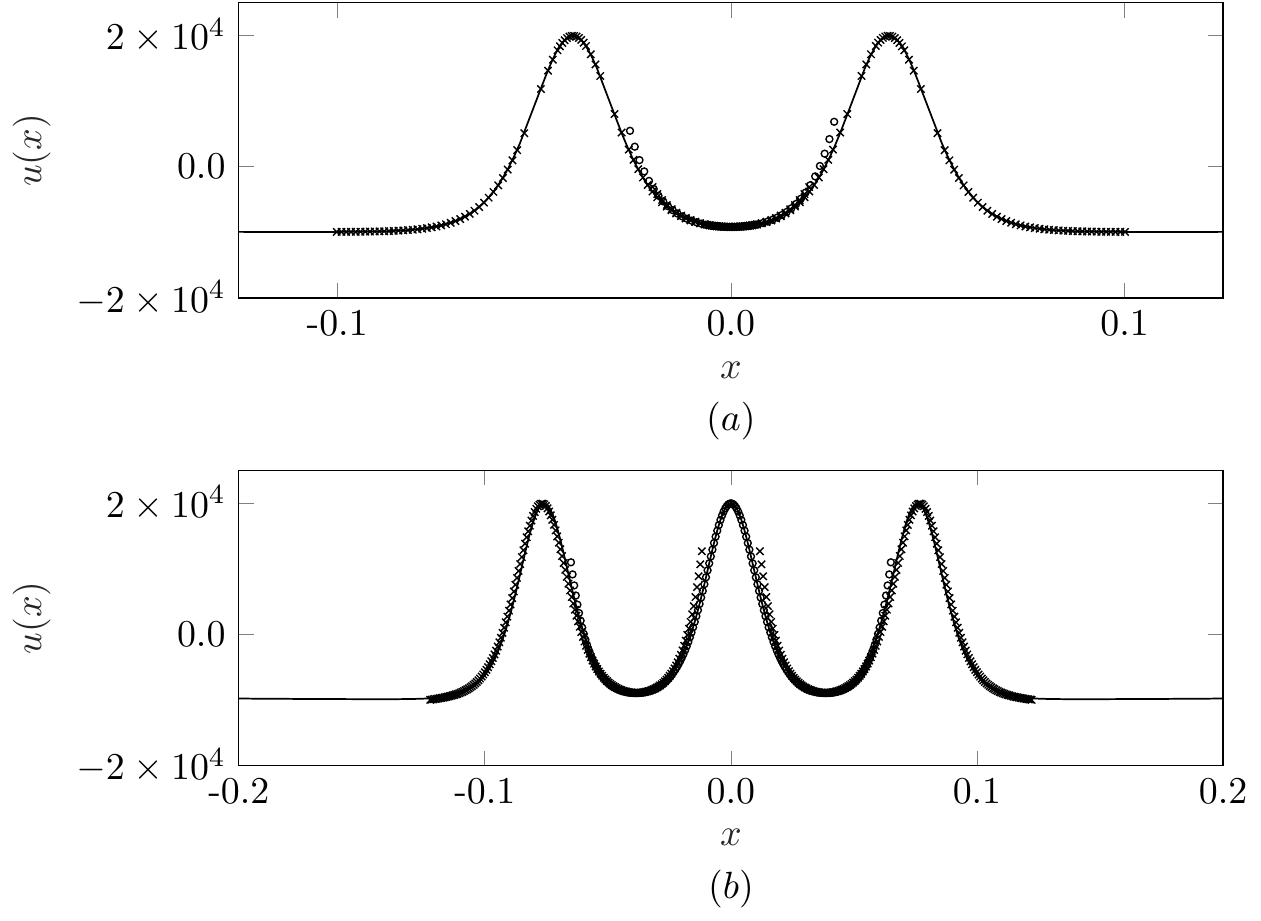}
\\
\caption{Double (upper panel) and triple (lower panel) homoclinic solutions for the Gaussian at $\alpha=10^8$. Comparison between the asymptotic predictions, shown with symbols,
and the numerical solution of \eqref{maineq}-\eqref{mainbc}, shown with a solid line. The circles correspond to $u=\alpha^{1/2}U_0 + U_1$ with $U_0$ given by $U_0^H$ or $-1$, and $U_1$ given by \eqref{shellington} or \eqref{gupJ} for the triple/double homoclinic respectively. The crosses correspond to $u=\alpha^{1/2}U_0 + \theta_1$ with 
 $U_0$ given by $U_0^H$ and $\theta_1$ given by \eqref{gupC} or \eqref{humma} for the triple/double homoclinic respectively. For the triple and double homoclinics $Y_1$ was taken to be the numerical solution of \eqref{manta} (with $\Lambda_1$ given by the equation in \eqref{manta2}) and \eqref{haggis}, respectively.}
\label{fig:clinic_comp}
\end{figure}


%

The double homoclinic corresponds to $\Lambda_1=0$ or
\bea \label{haggis}
\frac{12}{\delta^2}\ee^{-2\sq \ye_1} = -\frac{\sq}{8}f''(0)\ye_1
\eea
in which case 
\bea \label{gupL}
\hspace{-0.5in} \ye_1 \sim \tfrac{1}{\sq} \log(1/\delta) - \tfrac{1}{2\sq} \log\left ( \tfrac{1}{\sq} \log(1/\delta)  \right ) + \tfrac{1}{2\sq}\log\left(\frac{-96}{\sqrt{2}f''(0)}\right) + o(1).
\eea
Otherwise the analysis continues in a manner similar to that presented in section \ref{sec:oddhomoclinic}. In this case the sequence is established as
\bea \label{gupM}
\Lambda_{n+1} = \Lambda_n + \frac{\sq}{8}f''(0)\ye_n, \qquad 
\ye_{n+1} = \frac{1}{\sq}\log\left(\frac{12}{\Lambda_{n}\delta^2}\right) + \ye_n,
\eea
and $\hat y_n = \ye_{n+1}-\ye_{n} + \hat y_{n+1}$ for $n=1,2,\cdots$. 
The sequence for each successive $n$ corresponds to the value $\ye_1$ such 
that $\Lambda_n=0$ and has homoclinics at $y=\pm \ye_k$ for $k=1,\ldots,n$.
We find 
\bea \label{gupN}
\ye_n \sim \tfrac{(2n-1)}{\sq} \log(1/\delta)
\eea
Taking into consideration Remark 3 we may now check the validity of the expansion \eqref{gupI}. By taking the square root of \eqref{haggis}
it is clear that the $1/\delta^2$ terms in \eqref{gupJ} are in fact of $O(\lambda)$ where
\bea \label{gup_temp}
\lambda = \frac{1}{\delta}\log^{1/2}\left(1/\delta\right).
\eea
The expansion \eqref{gupI} should therefore be adjusted accordingly; this adjustment is permitted due to the linearity of the perturbation problems at each order of approximation.
%
\begin{table}
\begin{subtable}
\small
\centering
{
\renewcommand{\arraystretch}{1.8}
\begin{tabular}{ccc|ccccccc|ccccc}
\multicolumn{15}{c}{Gaussian triple homoclinic}\\
\hline
$\alpha$ &&& $u_0$ && $u_0^A$ && $\delta u$ &&& $Y_{H}$ && $Y_1$ && 
$\delta Y$\\
\hline
$10^5$ &&&  $599.92$ && $609.43$ && \cmmnt{9.51}$1.59\times 10^{-2}$ &&& $5.50$ && 5.37 && \cmmnt{1.43} $2.34\times 10^{-1}$ \\
$10^6$  &&&  $1963.09$ && $1972.37$ && \cmmnt{9.28}$4.73\times 10^{-3}$ &&& $6.15$
&& 6.10 && \cmmnt{1.27}$1.97\times 10^{-1}$ \\
$10^7$  &&&  $6283.33$ && $6292.32$ && \cmmnt{8.99}$1.43\times 10^{-3}$ &&& $6.86$
&& 6.83 && \cmmnt{1.16}$1.65\times 10^{-1}$ \\
$10^8$ &&&  $19954.49$ && $19963.16$ && \cmmnt{8.67}$4.34\times 10^{-4}$ &&& 7.58
&& $7.57$ && \cmmnt{1.07}$1.40\times 10^{-1}$\\
\\[-0.25in]
\multicolumn{15}{c}{Lorentzian triple homoclinic}\\
\hline
$\alpha$ &&& $u_0$ && $u_0^A$ && $\delta u$ &&& $Y_{H}$ && $Y_1$ && 
$\delta Y$\\
\hline
$10^5$ &&&  $602.66$ && $609.43$ && \cmmnt{6.77}$1.12\times 10^{-2}$ &&& $5.57$
&& 5.37 && \cmmnt{1.50}$2.34\times 10^{-1}$ \\
$10^6$ &&&  $1964.37$ && $1972.37$ && \cmmnt{8.00}$4.07\times 10^{-3}$ &&& $6.18$
&& 6.10 && \cmmnt{1.30}$1.96\times 10^{-1}$ \\
$10^7$ &&&  $6283.90$ && $6292.32$ && \cmmnt{8.42}$1.34\times 10^{-3}$ &&& $6.87$
&& 6.83 && \cmmnt{1.17}$1.65\times 10^{-1}$ \\
$10^8$ &&&  $19954.73$ && $19963.16$ && \cmmnt{8.43}$4.22\times 10^{-4}$ &&& 7.60
&& $7.57$ && \cmmnt{1.09}$1.40\times 10^{-1}$\\
\\[-0.25in]
\multicolumn{15}{c}{Gaussian double homoclinic}\\
\hline
$\alpha$ &&& $u_0$ && $u_0^A$ && $\delta u$ &&& $Y_{H}$ && $Y_1$ && 
$\delta Y$\\
\hline
$10^6$ &&&  $-755.32$ && $-764.44$ && \cmmnt{9.12}$1.21\times 10^{-2}$ &&& $3.29$
&& 3.27 && \cmmnt{0.02}$6.08\times 10^{-3}$ \\
$10^8$ &&&  $-9163.49$ && $-9174.92$ && \cmmnt{11.43}$1.25\times 10^{-3}$ &&& $4.01$
&& 4.01 && \cmmnt{0.00}$0.00$ 
\end{tabular}
\normalsize
}
\caption{Comparison between numerical and asymptotic results for the homoclinic glueing. In the table $u_0$ is the numerical estimate and $u_0^A$ is given by the pertinent asymptotic formula \eqref{asym_u0}. The location of the first homoclinic maximum is given numerically by $Y_H$ and the asymptotic estimate $Y_1$ is taken to be the numerical solution of \eqref{manta} (with $\Lambda_1$ given by the equation in \eqref{manta2}) for the triple homoclinic and as the numerical solution of \eqref{haggis} for the double homoclinic. The relative errors are $\delta u = |(u_0-u_0^A)/u_0|$ and ${\delta Y = |(Y_H-Y_1)/Y_H|}$.}
\label{table}
\end{subtable}
\end{table}
%

\vspace{0.125in}
The odd/even homoclinic analysis above is valid provided that the number of homoclinics $n$ is such that $\delta^2 Y_n\ll 1$. Given \eqref{bluh} and \eqref{gupN} this implies
\bea \label{gupP}
n \ll \frac{1}{\delta^2\log(1/\delta)}
\eea
as the condition that all the homoclinics are located where $|x|\ll 1$.
Since $f''(0)=-2$ for both the Gaussian and the Lorentzian, the results in the previous subsections give
the following approximations which may be applied to either case,
\begin{align}\label{asym_u0}
\nonumber
&u(0) \sim 2\alpha^{1/2}  - 4(\log 2)  \qquad   & (\mbox{single homoclinic}) \\
&u(0) \sim -\alpha^{1/2} + 2^{5/4}3^{1/2} \alpha^{1/4} \ye_1^{1/2} & (\mbox{double homoclinic}) \\
\nonumber
&u(0) \sim 2\alpha^{1/2}  - 4\log \alpha^{1/2} \qquad & (\mbox{triple homoclinic}) 
\end{align}
In Table \ref{table} we compare these asymptotic predictions with numerical calculations for the triple homoclinic for the Gaussian and the Lorentzian, and the double homoclinic for the Gaussian. In figure~\ref{fig:clinic_comp} we show a comparison between the asymptotic homoclinic glueing predictions and numerical solutions for the Gaussian forcing that demonstrate strong agreement between the two.

\section{Discussion}\label{sec:discuss}

We have analysed solutions to the problem \eqref{maineq}-\eqref{mainbc} for the case of a top hat forcing, a Gaussian forcing and a Lorentzian forcing, with particular attention paid to the limits of small and large $\alpha$. 
We have presented an asymptotic construction to provide supporting evidence for the existence of termination points on the solution branches for 
forcing functions which decay in the far-field faster than $1/x^4$, which includes the Gaussian forcing. We have also presented an asymptotic description 
of the large $\alpha$ solution profiles using the method of homoclinic glueing, which can be applied to any smooth forcing with a local 
maximum.
\begin{figure}[!t]
\centering
\subfigure[]{\includegraphics[width=3.0in]{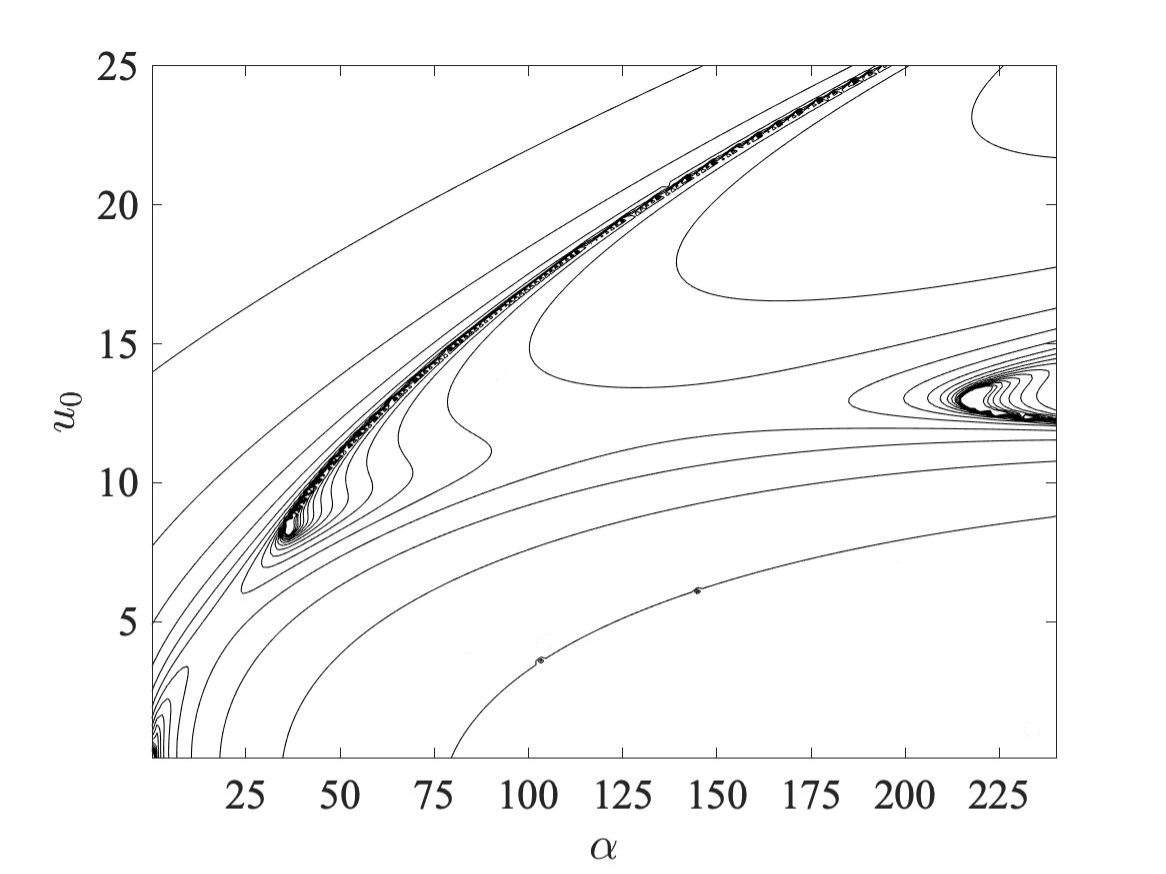}}
\subfigure[]{\includegraphics[width=3.0in]{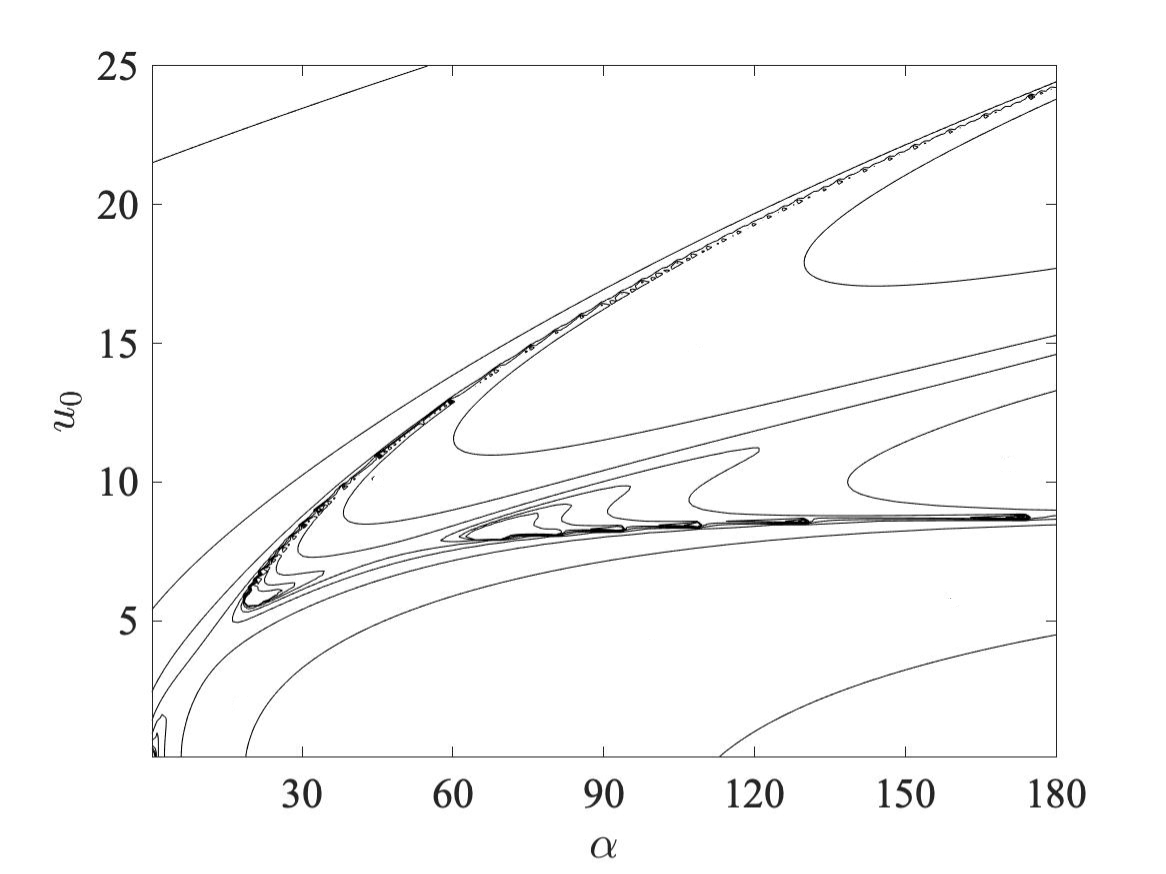}}
\subfigure[]{\includegraphics[width=3.0in]{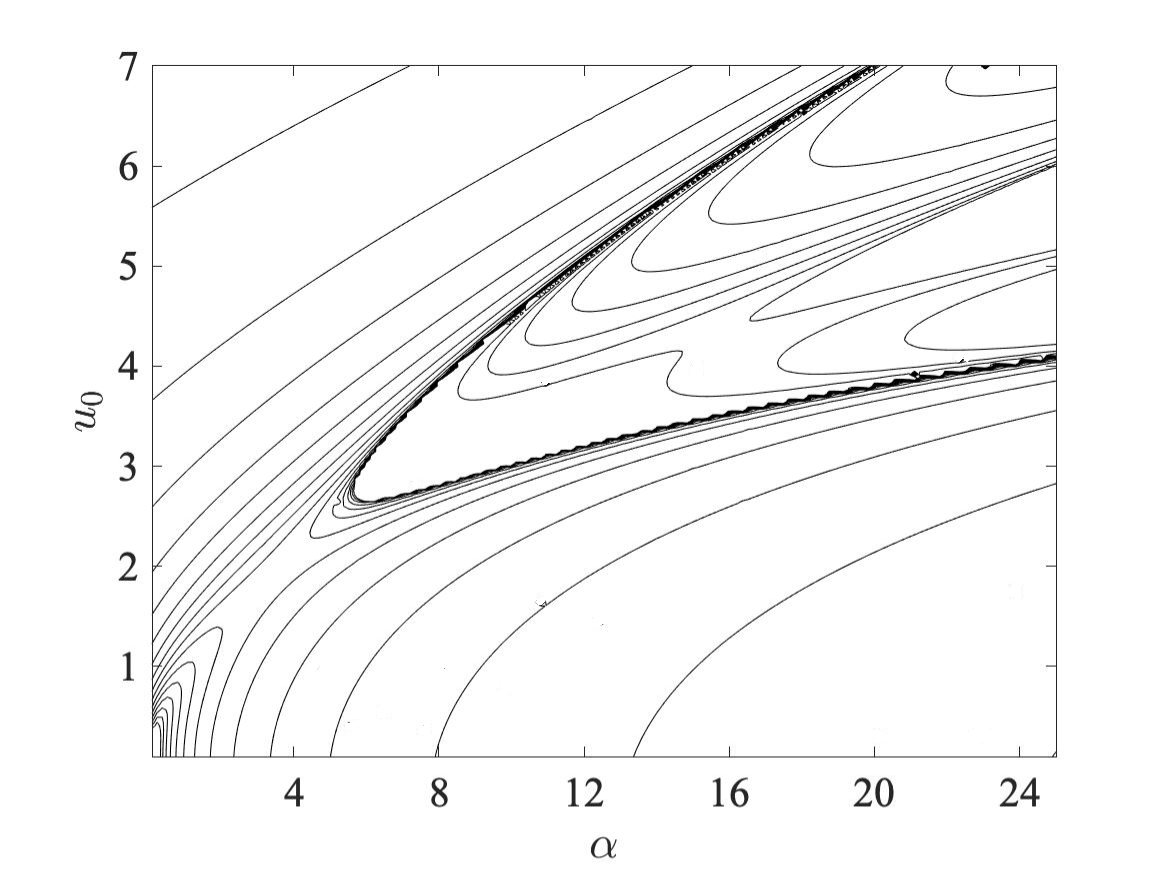}}
\caption{Contours of the blow-up point $x_0$ in the $(\alpha,u_0)$ plane for the hybrid forcing \eqref{synthesis} for (a)  $a=0$, (b) $a=0.7$ and (c) $a=1$.}
\label{fig:contours}
\end{figure}

The structure of the solution space is similar for all three forcing functions, each with an apparently infinite number of solution branches with qualitatively similar features in the solution profiles on each branch. Of particular note, however, is the presence of termination points for the Gaussian forcing on all branches, and the linking together of branches $B_n$, $B_{n+1}$ for the Lorentzian forcing (which decays more slowly that $1/x^4$ in the far-field). Some further insight into these different characteristic features can be obtained by attempting to continuously deform one forcing function into another. This can be achieved by considering the hybrid forcing function
\bea \label{synthesis}
f(x) = \frac{\ee^{-(1-a)x^2}}{1+ax^2}
\eea
for $0\leq a \leq 1$. As was noted above, the generic far-field behaviour for $u(x)$ in \eqref{maineq} corresponds to blow-up at the finite value $x=-x_0$ via \eqref{algdec}.
Figure \ref{fig:contours} shows contours of the blow-up point $x_0<0$ for the forcing \eqref{synthesis} at three sample values of $a$. Solution branches can be discerned on which formally $x_0=-\infty$ to satisfy the far-field condition \eqref{mainbc}. Of particular note when $a<1$ are the two saddle points in the contour map that are located for $a=0.7$ roughly at $(\alpha, u_0)=(11,3.6)$ and $(48,7.5)$. These saddle points persist on decreasing $a$ to zero, and indeed the contour plot for the Gaussian forcing is qualitatively similar to Figure \ref{fig:contours}. As $a$ is increased towards unity the tips of the two solution branches move toward each other. They pinch together at the rightmost of the two saddle points when $a=1$ to form the continuous solution branch labelled $B_1$, $B_2$ for the Lorentzian forcing in Figure \ref{fig:BVP_agnesi1}.

\begin{figure}
\centering
\includegraphics[width=4.0in]{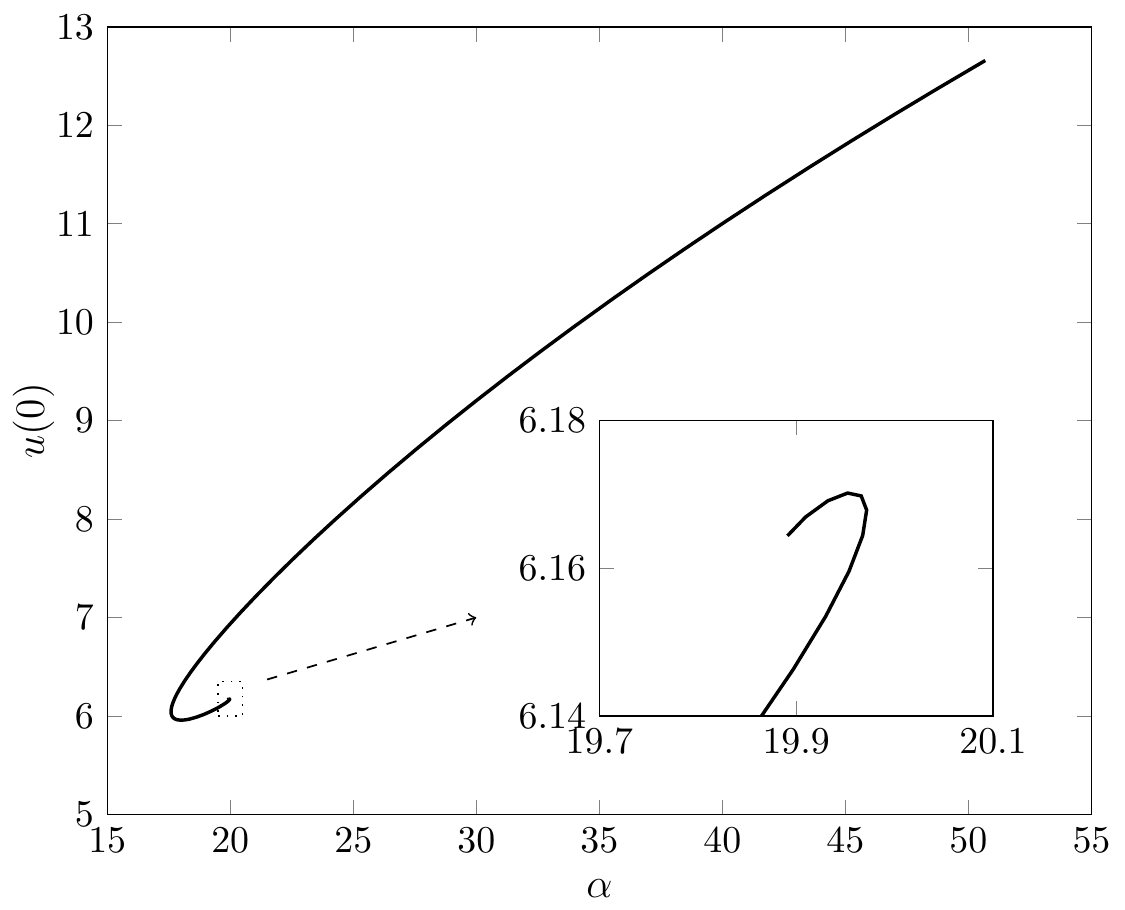}
\caption{Solution branch $B_1$ for the marginal case forcing \eqref{marginalf}.}
\label{fig:B1_marginal}
\end{figure}

As has been discussed, the existence of termination points at non-zero $\alpha$ appears to hinge on the decay rate of the forcing in the far-field relative to the inverse fourth power (see equation \ref{fork}). A forcing function for which $f(x)=O(1/x^4)$ as $x\to \infty$, for example,
\bea \label{marginalf}
f(x) = \frac{1}{1+x^4},
\eea
presents a marginal case. For this forcing, the large $x$ behaviour of the solution is
\bea
u \sim m_{\pm}/x^2, \qquad m_{\pm} = -3\pm (9+\alpha)^{1/2} \quad (m_+>0, \quad m_-<0),
\eea
constituting a balance between all three terms in (\ref{maineq}). The branch $B_1$ for the forcing \eqref{marginalf} is shown in Figure \ref{fig:B1_marginal}. The inset suggests that the branch terminates at $\alpha\approx 19.9$; in fact the branch spirals inwards toward the termination point beyond the last point reached by our numerics (see \ref{sec:marginal} for details). This underscores the subtle behaviour that can be found in problems of this type.

While we have used the particular class of equations \eqref{maineq} as an example to illustrate the idea of a termination point and to demonstrate the use of the method of homoclinic glueing, we believe that this class is simple enough to act as a paradigm for a much broader set of problems.
Finally we note that the homoclinic glueing analysis presented here is valid provided that the number of homoclinics satisfies condition \eqref{gupP}. If this condition is violated then a different approach is needed. This is the subject of ongoing work.


\ack
MGB and JRK gratefully acknowledge support from the ICMS 
as follow-on funding from the workshop `Applied and Computational Complex Analysis' held in Edinburgh from 8th-12th May 2017, the LMS under the Scheme 4 Grant 41853, and from the Isaac Newton Institute, Cambridge, as part of the programme `Complex analysis: techniques, applications and computations'  held from 2nd September 2019 to 19th December 2019. JSK gratefully acknowledges support from the IMA under the Small Grant Scheme.


\appendix

\section{Further details for the homoclinic glueing}\label{app:glue}

In this Appendix we provide some additional details of the calculations performed in the homoclinic glueing of section \ref{sec:glue}. 

\subsection*{The homoclinc glueing problem for $U_1$}

The solution of first order homoclinc glueing problem (\ref{kwazii}), namely
\bea \label{kwazii2}
\frac{\dd^2 U_1}{\dd y^2} + 2(3\mbox{sech}^2(y/\sqrt{2}) - 1)U_1 = \tfrac{1}{2}f''(0)y^2,
\eea
can be readily obtained readily by making the substitution $t=\tanh (y/\sq)$. The solution is found to be
\bea
U_1 = U_{10}\, \Phi_s(y/\sq) + \tilde U_{10}\, \Phi_a(y/\sq) + \tfrac{1}{2}f''(0)\Phi_p^{(2)}(y/\sq),
\eea
for arbitrary constants $U_{10}$ and $\tilde U_{10}$. (Note to fulfil \eqref{kwazii_cond} we set $\tilde U_{10}=0$.)
The antisymmetric and symmetric complementary functions are
\begin{align}
\hspace{-0.05in}
\Phi_a(y) = \tnh \, \ssq, \qquad  \Phi_s(y) =  \frac{1}{8}\Big( 15\ssq - 2\cosh^2 y - 5  - 15y\Phi_a(y)  \Big).
\end{align}
and they satisfy 
\bea
\Phi_a(0) = 0, \quad \frac{\dd \Phi_a}{\dd y}(0) = 1, \qquad 
\Phi_s(0) = 1, \quad \frac{\dd \Phi_s}{\dd y}(0) = 0.
\eea
The particular integral form is given by
\bea \label{carrot}
\hspace{-0.1in}
\Phi_p^{(2)}(y) = 2\Big [y \Big ( y\,\ssq  - 2\tanh y\Big ) - \log \ssq \Big ]\, \Phi_s(y)  + \Phi_a(y) I(y),
\eea
where
$$
I(y) = 4\int_{0}^{y} s^2\,\Phi_s(s) \,\dd s,
$$
and it satisfies 
\bea
\Phi_p^{(2)}(0)=\frac{\dd \Phi_p^{(2)}}{\dd y} (0) =0.
\eea
It will be helpful to note that 
\bea
4\int_{0}^{\infty} s^2\left(\Phi_s(s) + \tfrac{1}{16}\ee^{2s} + \tfrac{3}{4} \right )\dd s
= -\tfrac{1}{16} \left(5\pi^2 + 1\right). 
\eea

For reference we note the asymptotic properties of these functions. As $y\to +\infty$
\bea \label{app:asym}  
\hspace{-0.25in}
\Phi_a(y/\sq) \sim 4\ee^{-\sq y} + o\big(\ee^{-\sq y}\big), \qquad 
\Phi_s(y/\sq) \sim -\tfrac{1}{16}\ee^{\sq y} + o\left (\ee^{\sq y}\right), 
\eea
and
\bea \label{app:asym2}
\hspace{-0.25in}
\Phi_p^{(2)}(y/\sq) \sim \tfrac{1}{4}(\log 2)\, \ee^{\sq y} - \tfrac{1}{2}(y^2+1) + O(1);
\eea
and as $y\to -\infty$ we have
\bea \label{app:asym3}
\nonumber
&& \Phi_a(y/\sq) \sim -4\ee^{\sq y} + 16\,\ee^{2\sq y} + o\left (\ee^{2\sq y}\right), \\
&& \Phi_s(y/\sq) \sim -\tfrac{1}{16}\ee^{-\sq y} - \tfrac{3}{4} + \left(\tfrac{119}{16}+\tfrac{15}{2\sqrt{2}}y\right)\ee^{\sq y} + o\left (\ee^{\sq y}\right),\\
\nonumber
&& \Phi_p^{(2)}(y/\sq) \sim \tfrac{1}{4}(\log 2)\, \ee^{-\sq y} - \tfrac{1}{2}(y^2+1) + 3\log 2 + q(y/\sq) \ee^{\sq y},
\eea
where 
\bea
q(y) = 4y^3 - 5y^2 + \left(\tfrac{11}{2}-30\log 2\right)y + \tfrac{5\pi^2}{4} - \tfrac{19}{8} - \tfrac{119}{4}\log 2.
\eea

\subsection*{The homoclinic glueing problems for $\phi$ and $\psi_k$}

Herein we provide details of the first order problems satisfied by the functions $\phi(y_1)$ and $\psi_k(y_1)$ for $k=0,1,2$ which appear in the general solution for $\theta_1$ in (\ref{gupC}). The problem for $\phi$ is
\bea
\frac{\dd^2 \phi}{\dd y_1^2} + 2(3\mbox{sech}^2(y_1/\sqrt{2}) - 1) \phi = 0.
\eea
Using the results from above, the general solution may be written as
\bea
\phi = C_a \Phi_a(y_1/\sq) + C_s \Phi_s(y_1/\sq)
\eea
for constants $C_a$, $C_s$. Inspecting the asymptotic forms \eqref{app:asym} and \eqref{app:asym2} we see that to fulfil
the glueing condition \eqref{octopod} we must set $C_a = 0$ and $C_s=-16$ and so the solution for $\phi$ is even in $y_1$. 

The functions $\psi_k$ satisfy the problems 
\bea
\frac{\dd^2 \psi_k}{\dd y_1^2} + 2(3\mbox{sech}^2(y_1/\sqrt{2}) - 1)\psi_k = y_1^k
\eea
for $k=0,1,2$. The solutions are:
\begin{align} \label{toadstool}
\psi_k(y_1) = D_{a}^{(k)}\Phi_a(y_1/\sq) + D_{s}^{(k)} \Phi_s(y_1/\sq) + \Phi_p^{(k)}(y_1/\sq),
\end{align}
where the $D_{a}^{(k)}$ and $D_{s}^{(k)}$ are arbitrary constants and
\begin{align}
\hspace{0.in} \Phi_p^{(1)}(y_1) = \tfrac{\sqrt{2}}{8}\Big( \!  
\sinh (2y_1) -4y_1  + 6\tanh y_1 - (6y_1^2+15)\Phi_a(y_1) + 12y_1\,\sech^2 y_1 
\!\Big)
\end{align}
and
\begin{align}
\nonumber
\Phi_p^{(0)}(y_1) = \tfrac{1}{8}(15\tanh y_1 + 3y_1)\Phi_a(y_1) & - \tfrac{1}{4}\tanh^2\! y_1 \\
&- \tfrac{1}{8}\sinh^2\!y_1 (15\tanh^4 \!y_1 - 25\tanh^2 \!y_1 + 8).
\end{align}
We note the following asymptotic properties. As $y_1\to +\infty$
\bea
\hspace{-0.85in}
\Phi_p^{(1)}(y_1/\sq) \sim  \tfrac{\sqrt{2}}{16} \ee^{\sq y_1} - \tfrac{1}{2}y_1 + O(1), \qquad \Phi_p^{(0)}(y_1/\sq) \sim \tfrac{1}{16}\ee^{\sq y_1} + O(1);
\eea
and as $y_1\to -\infty$
\bea
\nonumber
& \hspace{-0.85in} \Phi_p^{(1)}(y_1/\sq) \sim -\tfrac{\sqrt{2}}{16} \ee^{-\sq y_1} - \tfrac{1}{2}y_1 
- \tfrac{3\sqrt{2}}{4} + \sqrt{2}\left(-\tfrac{95}{16}+9 y \right)\ee^{\sq y} + o\left(\ee^{\sq y}\right),\\
\\
\nonumber
& \Phi_p^{(0)}(y_1/\sq) \sim \tfrac{1}{16}\ee^{-\sq y_1} + \tfrac{1}{4} - \left(\tfrac{3}{2}y+\tfrac{23}{16}  \right)\ee^{\sq y} + o\left(\ee^{\sq y}\right),
\eea
Taking account of these asymptotic forms and those given above, the 
solutions that adhere to the glueing conditions \eqref{cbeebies} are given by 
\eqref{toadstool} with
\begin{align}
&D_{a}^{(0)} = -\tfrac{3}{2}, \qquad D_{a}^{(1)} = \tfrac{107\sqrt{2}}{8}, \qquad D_{a}^{(2)} = \tfrac{5\pi^2}{16} - \tfrac{19}{32}, \\
& D_{s}^{(0)} = 1, \qquad \quad \! D_{s}^{(1)} = -\sqrt{2}, \qquad \, D_{s}^{(2)} = 4\log 2.
\end{align}


\section{Stokes line analysis for the Lorentzian forcing}

The difficulty noted by Chapman {\it et al.} \cite{CKA}, for example, requires an analysis of the Stokes lines that emanate from the singularities of the leading order term in the expansion \eqref{gayle} extended into the complex plane. With this in mind, when referring to the analysis in section \ref{sec:agnesi_large} we shall replace $x$ with $z\in \mathbb{C}$.

The expansion \eqref{gayle} is a divergent asymptotic series whose terms
are generated by the recurrence relation \eqref{rashid}.
Following Chapman {\it et al.} \cite{CKA} we optimally truncate the asymptotic series at its smallest term, writing
\begin{equation}
W(z) = \sum_{n=0}^{N}\mu^n W_{n}(z) + R_{N}(z), \qquad W_n(z)=m^{2n+1}p_{2n}(z)(1+z^2)^{-(1+3n)/2},
 \label{kohli}
\end{equation}
where $R_N(z)$ is a remainder term and $m=\pm 1$ corresponding to the choice of sign made in \eqref{cricket}. The optimal truncation level $N$ follows from knowledge of the large $n$ behaviour of $W_n(z)$. We make the usual ansatz, writing (cf. Dingle \cite{dingle}),
\begin{equation}
W_n(z)\sim A(z)\,\frac{\Gamma(2n+\gamma+1)}{\sigma^{2n+\gamma+1}}
\label{joss}
\end{equation}
as $n\to \infty$, where $\Gamma(z)$ is the Gamma function and the functions $A(z)$, $\gamma(z)$ and the singulant
$\sigma(z)$ are all to be found. Substituting \eqref{joss} into
(\ref{rashid}) the balances at leading order, first order and second order determine that  (see Keeler \cite{keeler})
\bea \label{morgan}
{\sigma'}^2 = -2W_0(z)
\eea
where $'$ means $\mbox{d}/\mbox{d}z$, and that $\gamma(z)=-1/6$ and 
$A(z) = \Lambda/(\sigma')^{1/2}$, where $\Lambda$ is the Stokes multiplier that will be determined below. Since $W_0$ has singularities at $z=\pm \ri$, then $W_n(z)$ will also have singularities at these locations, for all $n$. We shall focus on the singularity at $z=\ri$, the analysis for $z=-\ri$ being similar. Integrating \eqref{morgan}, the singulant takes the form
\begin{equation}
    \sigma(z) = \left\{\begin{array}{rcl}
     -\sqrt{2}\int_{\ri}^{z}(1+p^2)^{-1/4}\,\mbox{d}p & \mbox{if} & m = -1, \\[0.1in]
     \mathrm{i}\sqrt{2}\int_{\ri}^{z}(1+p^2)^{-1/4}\,\mbox{d}p & \mbox{if} & m = 1, 
\end{array}\right.  
\label{buttler}
\end{equation}
where the lower integration limit has been chosen so that $\sigma(\ri)=0$. Finally, consistency of $W_n(z)$ for large $n$ between the forms given in \eqref{joss} and in \eqref{kohli} demands that (Keeler \cite{keeler})
\begin{equation}
\Lambda = \frac{2^{7/6}\pi^{1/2}}{3^{2/3}\Gamma\left(\frac{1}{3}\right)}\approx -0.7140572.
\end{equation}

\begin{figure}
\centering
\includegraphics[scale=0.9]{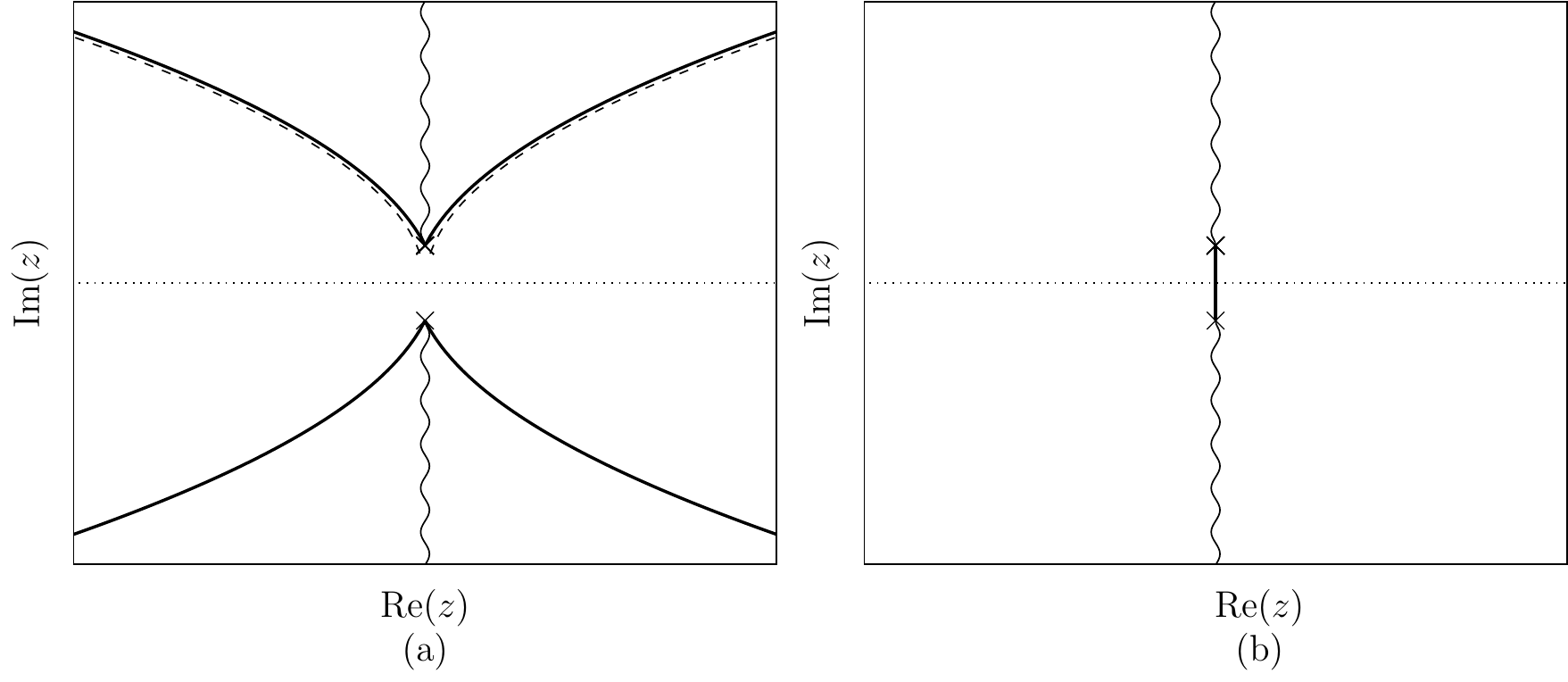}
\caption[Diagram showing Stokes' lines for `Witch of Agnesi'
Forcing]{(a) $m=-1$: numerically computed Stokes' lines $\mbox{Im}(\sigma)=0$, shown with solid lines, and the  asymptotic approximation \eqref{limerick}, shown with broken lines. (b) $m=1$: A Stokes line passes down the imaginary axis between $z=\pm \mathrm{i}$ crossing the real axis at $z=0$. In both (a) and (b) the branch cuts along the imaginary axis up from $z=\ri$ and down from $z=-\ri$ are shown with thin curvy lines.}
\label{fig:3_5_2}
\end{figure}

Stokes lines emerge from the points $z=\pm \ri$ where $\sigma$ vanishes (cf. Heading \cite{heading}) and hence, according to \eqref{joss}, the late form of $W_n$ is singular. According to Dingle \cite{dingle} the Stokes lines are traced by delineating the curves in the complex plane on which successive terms of the late order form \eqref{joss}, namely $W_n$ and $W_{n+1}$, have the same phase. Equivalently (see, e.g. \cite{heading}) on a Stokes line,
\begin{equation}
  \mbox{Im}(\sigma) = 0.
  \label{usman}
\end{equation}
The angle at which the lines emerge from the singularities is determined as follows. We note that
\bea
\hspace{-0.95in}
\sigma\sim \frac{2^{9/4}}{3}\mathrm{i}^{7/4}(z-\mathrm{i})^{3/4} \quad (\mbox{for}\,\,\, m=-1), \qquad 
\sigma\sim \frac{2^{9/4}}{3}\mathrm{i}^{11/4}(z-\mathrm{i})^{3/4} \quad (\mbox{for}\,\,\, m=1),
\label{agnew}
\eea
as $z\to\mathrm{i}$. Since the original problem posed on the real line is symmetric about $x=0$, it is natural to take the branch cuts that stem from the branch points at $z=\pm \ri$ to extend up the imaginary axis from $z=\ri$ and down the imaginary axis from $z=-\ri$ respectively. If local to $z=\ri$ we write $z-\mathrm{i}=R\mbox{e}^{i\psi}$ and $\sigma=r\mbox{e}^{\mathrm{i}\theta}$, then we should insist that
\bea \label{penguin}
-\frac{3\pi}{2}\leq\psi < \frac{\pi}{2}.
\eea
Considering first the case $m=-1$, (\ref{usman}) holds locally if 
$\psi = 4k\pi/3 - 7\pi/6$
for integer $k$. So there are two Stokes lines exiting $z=\mathrm{i}$ on which $\psi=-7\pi/6$ and $\psi=\pi/6$. (Similarly two Stokes lines exit $z=-\ri$ such that $\psi=-5\pi/6$ and $\psi=-\pi/6$). By appropriately deforming the contour of integration it can be shown that
for large $|x|$ the Stokes lines are approximated by
\bea \label{limerick}
y \sim \pm \rho |x|^{1/2}, \qquad 
\rho = \left(\tfrac{2}{\pi}\right)^{1/2}\Gamma^2(\tfrac{3}{4}) \approx 1.198.
\eea
The numerically computed Stokes lines in the upper half plane are shown in figure~\ref{fig:3_5_2}(a) together with the asymptotic approximation \eqref{limerick}. We note that the thickness of the Stokes layers about each Stokes line is of order $\mu^{1/4}|\sigma|^{1/2}$. Since $|\sigma|$ grows like $|z|^{1/2}$ for large $x$, the Stokes layer thickness grows as $x^{1/4}$, and hence the Stokes layer cannot impinge on the real line. It follows that for $m=-1$ the Stokes phenomenon can be ignored and the expansion \eqref{gayle} holds for all real $x$.

The situation is different for $m=1$. In this case the local form \eqref{agnew} leads to the conclusion that there is a 
Stokes' line along the imaginary axis between $z=\mathrm{i}$ and $z=-\ri$. An analysis similar to that in Chapman {\it et al.} \cite{CKA} (see Keeler \cite{keeler} for details) shows that the exponentially small remainder term
\bea \label{badterm2}
R_N(x)  \sim 2^{3/4}\pi \Lambda\, \mu^{-5/12}\, \ee^{-\rho (2/\mu)^{1/2}}\,(1+x^2)^{1/8}\, \cos \psi(x), 
\eea
where $\psi(x) = \left(2/\mu\right)^{1/2}\int_{0}^{x}(1+p^2)^{-1/4}\,\mbox{d}p$,
is activated on the real axis where the Stokes line crosses it at $x=0$. This term grows algebraically in $x$ so eventually the expansion \eqref{gayle} breaks down and there is therefore no solution of the original boundary value problem which is approximated by \eqref{gayle} for all $x$.


\section{Termination point analysis for the marginal case}\label{sec:marginal}

Our discussion for branch termination points hinged on the far-field decay behaviour \eqref{fork}. A forcing for which $f(x)=O(1/x^4)$ as $x\to \infty$, for example,
\bea \label{marginalf2}
f(x) = \frac{1}{1+x^4},
\eea
presents a marginal case. For large $x$,
\bea \label{usimmarg}
u \sim \frac{\varphi(\xi)}{x^2},
\eea
with $\xi=\log x$, where $\varphi$ satisfies the nonlinear equation
\bea
\frac{\dd^2 \varphi}{\dd \xi^2} - 5\frac{\dd \varphi}{\dd \xi} + 6\varphi + \varphi^2 = \alpha.
\eea
This has the two constant solutions
$\varphi = \varphi_{\pm}\equiv  -3 \pm (9+\alpha)^{1/2}$. These may be viewed as equilibria in the $(\varphi,\dd \varphi/\dd \xi)$ phase plane, wherein, assuming that 
$25-8\sqrt{9+\alpha}>0$, $(\varphi_+,0)$ is an unstable node
and $(\varphi_-,0)$ is a saddle node. So the far-field decay in \eqref{usimmarg} is such that,  as $\xi\to \infty$, $\varphi \propto \mbox{exp}(-p_{\pm}\xi )$ near to $\varphi_+$, and $\varphi \propto \mbox{exp}(-\lambda_{\pm} \xi)$ near to $\varphi_-$, where 
\bea
\hspace{-1.75cm}
p_{\pm} = 
\frac{1}{2}\left(-5 \pm (25-8\sqrt{9+\alpha})^{1/2}\right ), \qquad
\lambda_{\pm} = 
\frac{1}{2}\left(-5 \pm (25+8\sqrt{9+\alpha})^{1/2}\right ).
\eea
(If $25-8\sqrt{9+\alpha}<0$ then $\varphi_+$ is an unstable spiral.)
We conclude that two boundary conditions are required at $x=\infty$ to remove both of the eigenvectors at the unstable node, $\varphi_+$, and this leaves no degrees of freedom to satisfy the boundary condition at $x=0$. We therefore expect to find a solution in this case
only for special (possibly discrete) values of $\alpha$, labelled $\alpha^*$. 
%
%
Only one boundary condition is needed at $x=\infty$ to remove the unstable eigenvector at the saddle node, $\varphi_-$, leaving one degree of freedom to satisfy the condition at $x=0$. 

Working as in section \ref{sec:term_finite}, we perturb about the special solution $u_*(x)$ at 
$\alpha=\alpha^*$, writing
\bea
u \sim u_*(x) + \epsilon u_1(x), 
\eea
for small $\epsilon=\alpha-\alpha^*$. Substituting into \eqref{maineq} we find that the perturbation $u_1(x)$ satisfies
\bea \label{cagney}
\frac{\dd^2 u_1}{\dd x^2} + 2u_* u_1 = \frac{1}{1+x^4},
\eea
with $\dd u_1/\dd x(0) = 0$ and $u_1(x) \to 0$ as $x\to \infty$.
For large $x$, we have $u_*\sim \varphi_+(\xi)/x^2$ and the complementary functions for \eqref{cagney} are
$1/(x^{2}x^{p_{\pm}})$.
Numerical computations suggest that $\alpha^*$ is such that $p_{\pm}$ are a complex conjugate pair, 
so that for large $x$,
\bea \label{marg_non}
u \sim \frac{1}{x^2} \left(\varphi_{+} + \epsilon A^{+} x^{5/2}x^{\ri \tau}
+ \epsilon A^{-} x^{5/2}x^{-\ri \tau} \right),
\eea
for complex constants $A_{\pm}$ (with $A^+=\overline{A^-}$) and $\tau=(8\sqrt{9+\alpha}-25)^{1/2}$. A single relation is needed between the constants $A_{\pm}$ to have a boundary value problem for $u_1$. 

\begin{figure}
\centering
\includegraphics[scale=0.45]{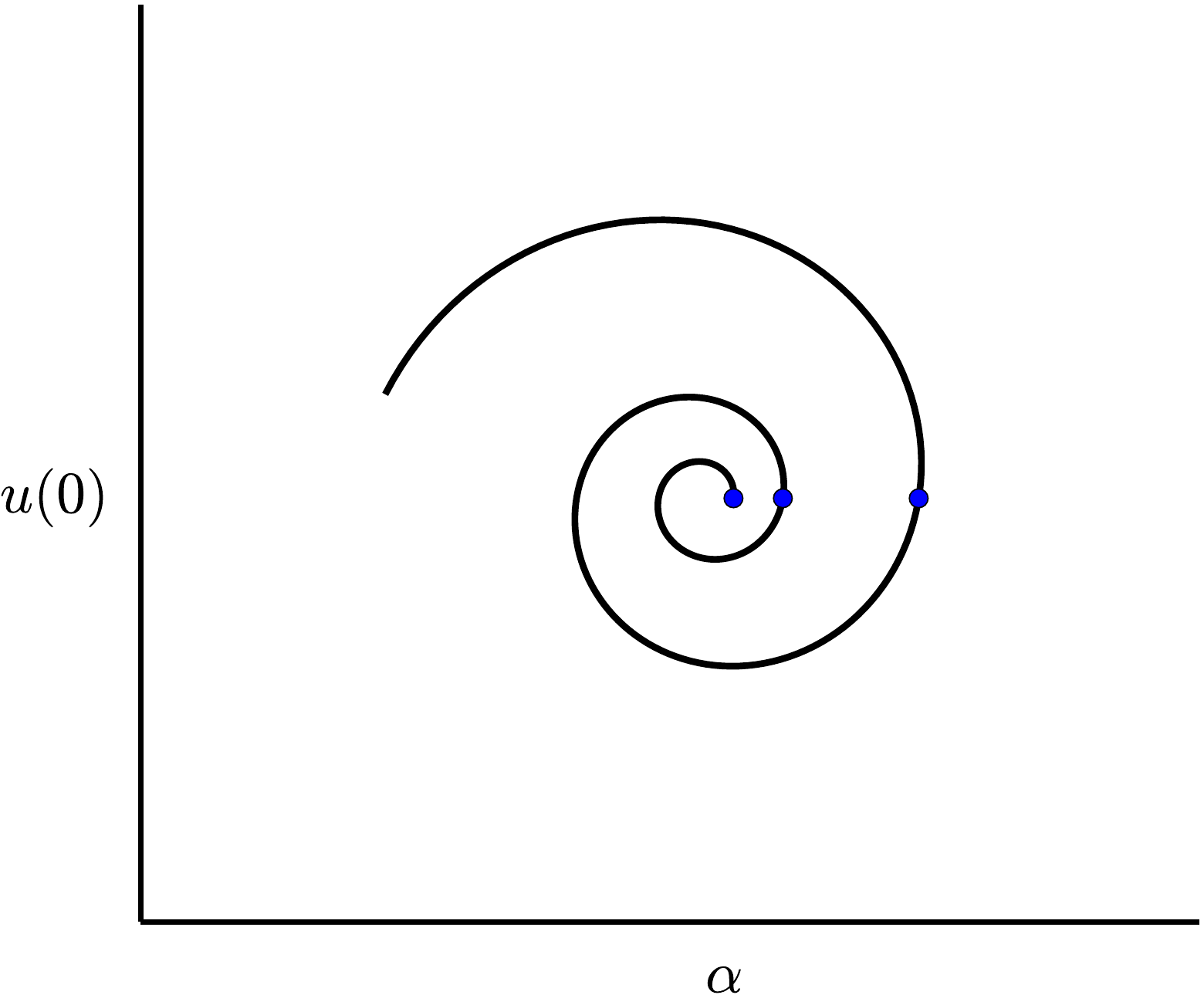}
\caption{Illustrative sketch of the local behaviour of the solution curve in the $(\alpha,u(0))$ plane near to $\alpha^*$ for the marginal case $f(x)=1/(1+x^4)$. The blue marker points indicate the change in the value $u(0)$ under the rescaling \eqref{epsrescale}.}
\label{fig:spiral}
\end{figure}
The non-uniformity of the expansion \eqref{marg_non} implies the presence of an outer region in which 
$x=\epsilon^{-2/5}X$ and $u=\epsilon^{4/5}U$ with $X=O(1)$ and $U=O(1)$. Writing $U = \varPhi/X^2$ and $\Theta=\log X$, we have
\bea \label{aaarggh}
\frac{\dd^2 \varPhi}{\dd \Theta^2} - 5\frac{\dd \varPhi}{\dd \Theta} + 6\varPhi + \varPhi^2 = \alpha^*
\eea
with
\bea \label{aaarggh2}
\varPhi \to \varphi_+ \quad \mbox{as $\Theta \to -\infty$},
\qquad
\varPhi \to \varphi_- \quad \mbox{as $\Theta \to \infty$}.
\eea
This will have a unique solution up to translations in $\Theta$ (compare travelling wave solutions to the Fisher-Kolmogorov equation, e.g. Britton \cite{britton})
with
\bea
\varPhi \sim \varphi_+ + B\ee^{5(\Theta+\Theta_0)/2}\cos\left(\tau[\Theta+\Theta_0]\right)
\eea
as $\Theta\to -\infty$ with $B$ real and effectively known from the unique solution to \eqref{aaarggh}, \eqref{aaarggh2}, and $\Theta_0$ arbitrary. Matching between the inner and the outer regions yields
\bea\label{sprout}
A_{\pm} = \frac{B}{2}\ee^{5\Theta_0/2}\exp\left (\mp \ri \tau \!\left[\Theta_0-\tfrac{2}{5}\log(1/\epsilon)\right] \right),
\eea
that is two equations which provide the relation between $A_+$ and $A_-$ alluded to above, and a condition to determine $\Theta_0$, on solving the boundary value problem for $u_1$.

According to \eqref{sprout} the rescaling
\bea \label{epsrescale}
\epsilon \mapsto \e^{-5\pi/\tau}\,\epsilon, \qquad \mbox{arg} A_{\pm} \mapsto \mbox{arg} A_{\pm} \pm 2\pi  
\eea
leaves the solution for $u_1$ unchanged. Therefore in this case the solution branch spirals into the termination point at $\alpha=\alpha^*$ in the manner sketched in figure \ref{fig:spiral}.



\end{document}